\documentclass{article}
\usepackage{amsmath}
\usepackage{amsfonts}
\usepackage{amsthm}
\usepackage{amssymb}
\usepackage[french]{babel}
\usepackage[latin1]{inputenc}
\usepackage[T1]{fontenc}
\usepackage[breaklinks=true]{hyperref}
\usepackage{eprint}

\newtheorem{thm}{Th\'eor\`eme}
\newtheorem{defi}{D\'efinition}
\newtheorem{cor}{Corollaire}

\newtheorem{lem}{Lemme}
\newtheorem{prop}{Proposition}
\newtheorem{rem}{Remarque}

\DeclareMathOperator{\car}{car}

\DeclareMathOperator{\MOD}{Qcoh}

\DeclareMathOperator{\rk}{rk}

\DeclareMathOperator{\AUT}{\underline{Aut}}

\DeclareMathOperator{\Hom}{Hom}
\DeclareMathOperator{\HOM}{\underline{Hom}}

\DeclareMathOperator{\Spec}{Spec}

\DeclareMathOperator{\Vect}{Vect}
\DeclareMathOperator{\VECT}{\mathcal{ VECT}}
\DeclareMathOperator{\PAR}{PAR}
\DeclareMathOperator{\Par}{Par}
\DeclareMathOperator{\coker}{coker}

\DeclareMathOperator{\obj}{obj}
\DeclareMathOperator{\Pic}{Pic}
\DeclareMathOperator{\Div}{div}
\DeclareMathOperator{\can}{can}
\DeclareMathOperator{\ch}{ch}
\DeclareMathOperator{\td}{td}
\DeclareMathOperator{\Coh}{Coh}
\DeclareMathOperator{\id}{id}
\DeclareMathOperator{\Func}{Func}
\DeclareMathOperator{\supp}{supp}
\DeclareMathOperator{\Sym}{Sym}
\DeclareMathOperator{\SPEC}{\bold{Spec}}

\input xy
\xyoption{all}

\begin{document}
 
\bibliographystyle{tocplain} 
\author{Niels Borne}
\title{Fibrés paraboliques et champ des racines \\(version préliminaire)}
 
 \maketitle

\section{Introduction}

\subsection{Résultats}

Soit $X$ une courbe projective et lisse sur un corps
algébriquement clos $k$, et $x$ un point fermé. Un
faisceau localement libre $\mathcal E$ sur $X$ est dit \emph{fini}
s'il existe deux polynômes distincts $P,Q$ à coefficients entiers positifs
tels que $P(\mathcal E)\simeq Q(\mathcal E)$, l'exemple classique
étant donné par les faisceaux inversibles sur $X$ dont une puissance
tensorielle est triviale. Nous nous conformerons à la tradition en
parlant de \emph{fibrés} finis. Nous nous intéressons aux liens bien
connus que ces fibrés ont avec le groupe fondamental de
$(X,x)$, ce par quoi nous entendrons toujours groupe
fondamental profini (\cite{SGA1}).

Dans \cite{NoriRFG}, Nori utilise la correspondance due à
A.Weil, lorsque $\car k =0$, entre fibrés finis sur $X$ et
représentations de $\pi_1(X,x)$, pour définir par
extension en caractéristique arbitraire, le \emph{schéma en groupe
fondamental de $(X,x)$}, comme groupe fondamental
$\pi(X,x)$ de la catégorie tannakienne des fibrés
essentiellement finis sur $X$, munie du foncteur fibre défini par
$x$. Ce nom est justifié par le fait que le groupe
proalgébrique $\pi(X,x)$ classifie les torseurs pointés
sur $X$, de groupe structurel un schéma en groupe fini sur $k$.

Dans un travail publié ultérieurement (\cite{NoriFGS}), Nori
montre que la notion de schéma en groupe fondamental a encore un
sens pour une courbe pointée $(X,D)$, ce groupe étant cette fois
défini comme le groupe fondamental de la catégorie tannakienne des
fibrés \emph{paraboliques} essentiellement finis sur $X$. Il utilise
ici sa propre définition des fibrés paraboliques, la décrivant comme
"une légère modification de celle de Seshadri". Il démontre, en suivant 
l'analogie avec le cas des courbes complètes, que ce groupe
fondamental classifie de même les torseurs pointés sur $X-D$, de
groupe structurel un schéma en groupe fini sur $k$.

Ce résultat a de l'intérêt même en caractéristique zéro,
puisqu'il donne une réinterprétation algébrique du groupe
fondamental étale de $X-D$, comme groupe de Tannaka de la catégorie
des fibrés paraboliques finis sur $(X,D)$. Ceci permet de reformuler
le problème classique (évoqué dans \cite{SGA1}) d'une preuve
"algébrique" (disons, sans utiliser de théorème de type GAGA) des
théorèmes de structure des groupes fondamentaux des courbes
algébriques. Par exemple, on peut essayer de comprendre
l'équivalence entre les fibrés paraboliques finis sur $({\mathbb
P}^1,\{0,1,\infty\})$ et représentations de $\widehat{F_2}$, le
complété profini du groupe libre à deux générateurs, même si ce
problème reste probablement très difficile.

Aussi notre attention s'est-elle concentrée sur deux problèmes : déterminer si le théorème de Nori sur les fibrés paraboliques se généralise en dimension supérieure, au moins en caractéristique zéro, et éclaircir le lien entre sa définition des fibrés paraboliques et celle, classique, due à Seshadri. Nous ne résolvons que le second problème dans ce travail (nous comptons revenir sur le premier dans un prochain article).

Nous montrons essentiellement deux résultats, qui sont valables en dimension quelconque. Donnés un entier $r\geq 1$, un schéma (noethérien) $X$ muni d'un diviseur de Cartier effectif $D$, le premier énoncé (théorème \ref{thmprinc}) 
stipule qu'il y a une équivalence de catégories tensorielles entre fibrés paraboliques (au sens de Seshadri) à poids multiples de $\frac{1}{r}$ et faisceaux localement libres sur un certain champ de Deligne-Mumford $\mathfrak X$ sur $X$, qu'on peut appeler le champ des racines $r$-ièmes de $D$, où, comme on l'imagine, la structure de champ est concentrée le long du support de $D$. 
Le second résultat (théorème \ref{thmdeg}) affirme que s'il l'on suppose de plus $X$ et $D$ projectifs et lisses 
sur un corps algébriquement clos de caractéristique première à $r$, 
la notion de degré parabolique (définie dans ce cadre de généralité par Maruyama et Yokogawa, \cite{MY}) coïncide avec la notion usuelle de degré pour le faisceau correspondant sur $\mathfrak X$.

Ces deux résultats permettent de ramener les questions concernant les fibrés paraboliques à des questions concernant des fibrés usuels. Dans le paragraphe \ref{courb}, nous montrons comment les notions usuelles sur les fibrés paraboliques sur une courbe marquée $(X,D)$ s'interprètent naturellement lorsque l'on considère le(s) champ(s) associé(s). Il s'ensuit par exemple que les fibrés paraboliques finis sont semi-stables de degré $0$, la preuve pour les fibrés usuels sur une courbe complète s'adaptant au cas d'une orbicourbe complète. Ceci permet d'obtenir des informations sur la structure des fibrés paraboliques finis sur la droite projective (théorème \ref{strucfibparfindrproj}). Ce résultat est en fait une conséquence du théorème de Nori, mais nous n'en faisons pas vraiment usage ici.

\subsection{Origines}

La définition de Nori (voir \cite{NoriFGS}) d'un fibré parabolique sur $(X,D)$ est, grossièrement, la donnée d'un fibré usuel sur $X-D$ auquel on recolle des fibrés définis sur un revêtement ramifié d'un voisinage des points du support de $D$. Dans l'esprit, c'est une définition très proche de celle des fibrés sur le champ $\mathfrak X$, qui est en quelque sorte le revêtement ramifié minimal de $X$ sur lequel on peut voir tous les fibrés paraboliques à poids multiples de $\frac{1}{r}$ sur $X$ comme des fibrés usuels. 

L'idée d'associer à un revêtement galoisien ramifié de courbes $Y\rightarrow X$, de groupe $G$, et à un $G$-faisceau cohérent sur $Y$, un faisceau parabolique sur $(X,D)$, où $D$ est le diviseur de branchement du revêtement, semble due à Biswas (\cite{Biswasorb}). C'est la procédure que nous adoptons, pour le ``revêtement ramifié'' $\mathfrak X\rightarrow X$,  pour associer à un fibré sur $\mathfrak X$ un fibré parabolique sur $(X,D)$. Biswas décrit également une procédure en sens inverse, et il est mentionné dans \cite{BBN} qu'il en résulte une équivalence tensorielle entre certains $G$-faisceaux sur $Y$ et les faisceaux paraboliques sur $(X,D)$. Nous nous distinguons ici de (\cite{Biswasorb}) en donnant, dans le cadre plus général des champs de racines, une expression plus explicite et intrinsèque de l'équivalence réciproque, en termes de cofins (``coends''). Cette expression rend d'ailleurs le caractère tensoriel évident.

Le champ des racines $r$-ièmes d'un faisceau inversible munis d'une section sur $X$ de Cadman et Vistoli (\cite{Cadman}, \cite{AGV}) est un outil essentiel. Tout aussi indispensables nous sont les revêtements cycliques uniformes de Arsie-Vistoli (\cite{AV}), dont les champs quotient sont les modèles locaux du champ des racines $r$-ièmes (voir théorème \ref{root}).

La définition des fibrés paraboliques que nous adoptons ici est due à Maruyama-Yokogawa (\cite{MY}, \cite{Yokogawa}), c'est une adaptation de celle de Seshadri sur les courbes (\cite{Seshadri}, voir aussi appendice \ref{defSeshadri}) à la dimension quelconque.

Notre preuve de l'égalité du degré d'un faisceau parabolique avec le degré du faisceau associé sur le champ des racines $r$-ièmes repose de manière essentielle sur la formule de Grothendieck-Riemann-Roch pour les champs de Deligne-Mumford due à Toën (\cite{ToenThese},\cite{Toen}). Nous regroupons les résultats que nous utilisons dans un long appendice (\ref{GRR}).

Enfin c'est une note de Edixhoven (\cite{Edix}) qui a éveillé notre attention sur la structure des fibrés paraboliques finis sur la droite projective.  

\subsection{Limites de cet article, et possibles applications}
\label{sec:limites}

Le théorème \ref{thmprinc} permet de réinterpréter et de compléter certains travaux antérieurs sur les fibrés paraboliques dans le cadre plus général des champs de Deligne-Mumford. Voici deux exemples. La définition des fibrés principaux à structure parabolique de groupe structurel donné de \cite{BBN} s'interprète comme celle des fibrés principaux de même groupe structurel sur un champ de racines (et cette interprétation nous semble plus naturelle que la définition originelle en termes de foncteurs tensoriels). Second exemple : on peut utiliser les différentes théories des anneaux de Chow pour les champs de Deligne-Mumford (par exemple celle de \cite{VistInv}) pour donner une définition des classes de Chern des fibrés paraboliques également plus naturelle que celle de \cite{BisChern}. Mentionnons à ce sujet que le ``caractère de Chern à coefficient dans les représentations'' de \cite{Toen} permet de calculer la $K$-théorie de la catégorie des fibrés paraboliques (dont les poids ont des dénominateurs bornés) en termes d'anneaux de Chow d'un champ de Deligne-Mumford (le champ d'inertie du champ des racines), ce qui semble hors de portée des seules techniques paraboliques. 

Toutefois, nous ne prenons pas le soin de décrire ces applications en détail, à la fois par manque de place, par volonté de se concentrer sur le lien entre fibrés paraboliques et groupe fondamental, et surtout parce qu'il apparaît a posteriori que la définition des fibrés paraboliques adoptée ici (celle de \cite{Yokogawa}) ne soit correcte que dans le cas d'un diviseur lisse (et il nous semble à présent que l'équivalence tensorielle de \cite{Biswasorb}, \cite{BBN} n'est valable que dans ce cadre).  En d'autres termes, bien que le théorème \ref{thmprinc} soit vrai pour un diviseur quelconque, elle n'est intéressante que pour un diviseur lisse. Nous comptons revenir sur le problème de la définition des fibrés paraboliques relatifs à un diviseur à croisements normaux dans un article ultérieur.
Il serait aussi très souhaitable de comprendre le théorème \ref{thmprinc} non 
seulement comme une équivalence de catégories, mais comme un équivalence de
champs, puisque cela permettrait de voir \cite{Yokogawa} sous un jour
nouveau. 

Mentionnons enfin le lien avec les structures logarithmiques et en particulier avec le travail de Matsuki et Olsson \cite{MatOls}, qu'il faudrait aussi expliciter.

\subsection{Remerciements}

Je tiens à remercier les nombreux mathématiciens  qui m'ont permis de mener à
bien cet article, et parmi eux, spécialement  I.Biswas et D. S.Nagaraj, pour
leurs patientes explications sur les fibrés paraboliques, B.Toën, pour de
multiples éclaircissements, A.Chiodo, pour d'intéressantes discussions sur les
champs de racines, et enfin, de manière générale, M.Emsalem et A.Vistoli.

 \section{Fibrés paraboliques  : la définition de Maruyama-Yokogawa}

On fixe un schéma noethérien $X$, et un diviseur de Cartier effectif $D$ sur $X$, et $r\geq 1$ un entier.

\subsection{Poids à dénominateurs multiples d'un entier donné}
On identifie tout ensemble ordonné à la catégorie correspondante.

On appelle $\MOD(X)$ la catégorie des faisceaux quasi-cohérents sur $X$.

Soit $\mathcal E_{\cdot} :
(\frac{1}{r}\mathbb Z)^{op}\rightarrow \MOD(X)$ un foncteur.
Pour  $i\in \frac{1}{r}\mathbb Z$
on dispose alors du décalage (shift) $\mathcal E[i]_\cdot$ défini sur
les objets par la formule
usuelle   $\mathcal E[i]_{j}=\mathcal E_{i+j}$, et pour $i\geq j$
de la transformation naturelle $\mathcal E[i]_\cdot\rightarrow
\mathcal E[j]_\cdot$.

\begin{defi}[\cite{Yokogawa}]
  \label{defYok}
  On définit la catégorie $\PAR_{\frac{1}{r}}(X,D)$ des fibrés paraboliques de
  poids multiples de $\frac{1}{r}$ comme ayant pour objets les
  couples  $(\mathcal E_{\cdot},j)$, où $\mathcal E_{\cdot} :
(\frac{1}{r}\mathbb Z)^{op}\rightarrow \MOD(X)$ est un foncteur, et
  $j:\mathcal E_{\cdot}\otimes_{\mathcal O_X}\mathcal
O_X(-D)\simeq \mathcal E[1]_\cdot$ est une transformation naturelle
  faisant commuter le diagramme :

 \begin{center}
\xymatrix@R=2pt{
  &\mathcal E_{\cdot}\otimes_{\mathcal O_X}\mathcal
O_X(-D)\ar[rr]^j\ar[dr]&& \mathcal E[1]_\cdot\ar[dl]& \\
  &&\mathcal E_{\cdot}&&}
\end{center}

Les morphismes de $(\mathcal E_{\cdot},j)$ à $(\mathcal
  E'_{\cdot},j')$ sont les transformations naturelles
  $\alpha :\mathcal E_{\cdot}\rightarrow \mathcal E'_{\cdot}$ telles
  que le diagramme suivant commute :

\xymatrix@R=2pt{
 && &\mathcal E_{\cdot}\otimes_{\mathcal O_X}\mathcal
O_X(-D)\ar[r]^{j}\ar[ddd]_{\alpha\otimes 1} & \mathcal
  E[1]_\cdot\ar[ddd]^{\alpha[1]} & \\
&&&&&\\&&&&&\\
  && &\mathcal E'_{\cdot}\otimes_{\mathcal O_X}\mathcal
O_X(-D)\ar[r]_{j'} & \mathcal E'[1]_\cdot&}

 \end{defi}

Pour alléger les notations, on omettra souvent de mentionner $j$ pour
désigner les objets de  $\PAR_{\frac{1}{r}}(X,D)$.

Lorsque $r=1$, il est clair que le foncteur d'oubli $\mathcal E_\cdot\rightarrow \mathcal E_0$ (ou évaluation en
zéro) $\PAR_1(X,D)\rightarrow \MOD(X)$ est une équivalence, dont une
équivalence réciproque est donnée par $\mathcal F \rightarrow \mathcal F
\otimes_{\mathcal O_X}\mathcal O_X(-\cdot D)$.

La catégorie  $\PAR_{\frac{1}{r}}(X,D)$ est munie de $\Hom$ internes.
Soient en effet $\mathcal E_{\cdot}$ et $\mathcal E'_{\cdot}$ deux
objets, on dispose déjà du faisceau
$\HOM(\mathcal E_{\cdot},\mathcal E'_{\cdot})$ défini par, pour tout ouvert $U$
de $X$ :
$\HOM(\mathcal E_{\cdot},\mathcal E'_{\cdot})(U)=\Hom(\mathcal
{E_{\cdot}}_{|U},{\mathcal E'_{\cdot}}_{|U})$.
On définit le faisceau parabolique $\HOM(\mathcal
E_{\cdot},\mathcal E'_{\cdot})_{\cdot}$ en posant pour tout $i\in
\frac{1}{r}\mathbb Z$ :

$$\HOM(\mathcal
E_{\cdot},\mathcal E'_{\cdot})_i= \HOM(\mathcal E_{\cdot},\mathcal E'_{\cdot}[i])$$

On définit alors le produit tensoriel $\mathcal
E_{\cdot}\otimes_{\mathcal O_X}\mathcal E'_{\cdot}$
de deux fibrés paraboliques $\mathcal E_{\cdot}$ et
$\mathcal E'_{\cdot}$ par la formule habituelle :

$$\HOM(\mathcal
E_{\cdot}\otimes_{\mathcal O_X}\mathcal E'_{\cdot}, \mathcal
E''_{\cdot})_ {\cdot}
\simeq \HOM(\mathcal E_{\cdot},\HOM(\mathcal E'_{\cdot},\mathcal E''_{\cdot})_{\cdot})_{\cdot}$$

On désigne par $\Vect(X)$ la sous-catégorie pleine de $\MOD(X)$ dont
les objets sont les faisceaux localement libres (de rang fini).

\begin{defi}[\cite{Yokogawa}]
  La catégorie $\Par_{\frac{1}{r}}(X,D)$ des $r$-fibrés paraboliques
 localement libres (de rang fini) est la sous-catégorie pleine de
 $\PAR_{\frac{1}{r}}(X,D)$ dont les objets sont les couples
 $(\mathcal E_{\cdot},j)$ vérifiant :

 \begin{enumerate}
   \item  $\mathcal E_{\cdot} : (\frac{1}{r}\mathbb Z)^{op}\rightarrow
   \MOD(X)$ se factorise en $(\frac{1}{r}\mathbb Z)^{op}\rightarrow \Vect(X)$.
   \item Pour $i\leq i' < i+1 $, $\coker(\mathcal E_{i'} \rightarrow \mathcal E_i)$
  est localement libre comme $\mathcal O_D$-module.
 \end{enumerate}

\end{defi}
  \begin{rem}[\cite{Yokogawa}]
\label{remfibrpar}  Donnés $(\mathcal E_{\cdot},j)$ et $(\mathcal E'_{\cdot},j')$ deux
  objets de $\Par_{\frac{1}{r}}(X,D)$,
  la compatibilité des transformations naturelles de
  $\mathcal E_{\cdot}$ vers $\mathcal E'_{\cdot}$ avec $j$ et $j'$ est
  automatique.
\end{rem}

\subsection{Poids rationnels quelconques}
Si $r|r'$, l'inclusion $(\frac{1}{r}\mathbb
Z)\subset(\frac{1}{r'}\mathbb Z)$ admet pour adjoint à gauche
le foncteur défini sur les objets par :
$\frac{l'}{r'}\rightarrow -\frac{1}{r}[-\frac{rl'}{r'}]$.

Ce dernier induit un foncteur fidèlement plein
$\PAR_{\frac{1}{r}}(X,D)\rightarrow \PAR_{\frac{1}{r'}}(X,D)$, adjoint
à gauche du foncteur de restriction $\PAR_{\frac{1}{r'}}(X,D)\rightarrow
\PAR_{\frac{1}{r }}(X,D)$.  En particulier pour $r=1$ on obtient un
foncteur $\MOD (X)\rightarrow \PAR_{\frac{1}{r'}}(X,D)$, adjoint à
gauche du foncteur d'oubli (ou évaluation en zéro)
$\PAR_{\frac{1}{r'}}(X,D) \rightarrow \MOD (X)$. Ceci permet de poser :

\begin{defi}
On définit
  la catégorie $\PAR(X,D)$ des fibrés paraboliques à
  poids rationnels comme

 $$\PAR(X,D)=\lim_{\stackrel{\longrightarrow}{r\in\mathbb N}} \PAR_{\frac{1}{r}}(X,D)$$

\end{defi}

\begin{rem}
  On peut aussi définir directement les objets de $\PAR(X,D)$ comme des
  couples  $(\mathcal E_{\cdot},j)$, où $\mathcal E_{\cdot} :
(\mathbb Q)^{op}\rightarrow \MOD(X)$ est un foncteur, et
  $j:\mathcal E_{\cdot}\otimes_{\mathcal O_X}\mathcal
O_X(-D)\simeq \mathcal E[1]_\cdot$ une transformation naturelle
  vérifiant la même relation de commutativité, avec la condition
  supplémentaire que la longueur de la filtration $(\mathcal
  E_{\alpha})_{\alpha \in [0,1]}$ est finie.
\end{rem}

On munit $\PAR(X,D)$ des Homs internes et du produit tensoriel induit
par ceux des $\PAR_{\frac{1}{r}}(X,D)$, qui sont évidemment
compatibles entre eux.

On définit la catégorie des fibrés paraboliques
 localement libres (de rang fini) sur $(X,D)$ de la manière évidente,
 à savoir :  $$\Par(X,D)=\lim_{\stackrel{\longrightarrow}{r\in\mathbb N}} \Par_{\frac{1}{r}}(X,D)$$

Comme on l'a vu, le foncteur d'oubli (ou évaluation en zéro)
$\Par(X,D)\rightarrow \Vect(X)$ a un adjoint à gauche, qu'on explicite :

\begin{defi}
Pour $\mathcal F\in \obj \Vect (X)$, on dispose du fibré parabolique
localement libre sur $(X,D)$, dit \emph{à structure spéciale}, et noté
$\underline{\mathcal F}.$, défini sur les objets par :
$$\underline{\mathcal F}_\cdot=   \mathcal F
\otimes_{\mathcal O_X}\mathcal O_X([-\cdot] D)$$
\end{defi}

L'utilité des structures spéciales vient en partie du fait classique suivant.

\begin{thm}

  \label{LocStrPar}
  Soit $\mathcal E_\cdot \in \obj \Par(X,D)$. Pour tout point $x\in
  X$, il existe un voisinage ouvert (pour la topologie de Zariski) 
$U$ de $X$, il existe une famille
  finie de nombres rationnels $(\alpha_a)_{a\in A}$ telle que

  $$(\mathcal E_\cdot)_{|U}\simeq \oplus_{a\in A}
  \underline{\mathcal O_U}_\cdot [\alpha_a]$$
\end{thm}

\begin{proof}
La démonstration est très proche de celle de la proposition \ref{StrucLoc}, que nous détaillons par la suite.
\end{proof}

 \subsection{Fibrés paraboliques et revêtements}

\begin{defi}[\cite{AV}]\label{RCU}
Soit $X$ un schéma noethérien. Un
\emph{revêtement cyclique uniforme de degré $r$} de $X$ est
un morphisme de schémas $\pi \colon Y \rightarrow X$, et une action du
schéma en groupes $\boldsymbol{\mu}_r$ sur $Y$, tels que pour tout
point $x$ de $X$, il existe un voisinage affine $U =  \Spec R$
de  $x$ dans $X$, un élément $s \in R$ qui n'est pas un diviseur de
zéro, et un isomorphisme de $U$-schémas $\pi^{-1}(U) \simeq
\Spec R[t]/(t^r - s)$ qui est $\boldsymbol{\mu}_r$-équivariant, où le
membre de droite est équipé de l'action évidente.
\end{defi}

L'ensemble de tels couples $(U,s)$ (resp. $(\pi^{-1}(U),t)$),
pour $U$ variant sur un recouvrement de $X$, définit un diviseur de
Cartier effectif $D$ sur $X$ (resp. $E$ sur $Y$),
que l'on appelle \emph{le diviseur de branchement}
(resp. \emph{le diviseur de ramification}) de $\pi$. Ils vérifient la
relation $\pi^*D=rE$.

\begin{rem}[\cite{AV}]\label{remRCU}
Donné $\pi \colon Y \rightarrow X$ un revêtement cyclique uniforme de
degré $r$, $\pi_*\mathcal O_Y$ est muni d'une $\mathbb Z/r$-graduation
naturelle, qu'on écrit
$\pi_*\mathcal O_Y =\mathcal L_0\oplus \mathcal L_1\oplus \cdots
\mathcal L_{r-1}$. En regardant localement on voit que les $\mathcal
L_j$ sont des faisceaux inversibles, et que le morphisme
canonique $\mathcal L_1^{\otimes r}\rightarrow \mathcal L_0
\simeq  \mathcal O_X$ est un
monomorphisme d'image $\mathcal O_X(-D)$.
Ainsi le faisceau $\mathcal L= \mathcal O_X(D)$ est-il muni d'une
racine $r$-ième canonique.

Réciproquement la donnée d'un faisceau inversible  $\mathcal L_1$ sur
$X$ et d'un monomorphisme  $\mathcal L_1^{\otimes r}\rightarrow
\mathcal O_X$ permet de définir un revêtement cyclique uniforme de
degré $r$ sur $X$, en considérant le schéma affine sur $X$
associé au faisceau d'algèbres $\frac{\Sym \mathcal L_1}{\mathcal L_1^{\otimes r}\rightarrow \mathcal O_X}$.

\end{rem}

\begin{defi}
  Soit $S$ un schéma de base.
  Soit $Y$ un $S$-schéma muni d'une action d'un $S$-schéma en groupes
  $G$. Un \emph{$G$-faisceau localement libre sur $Y$} est un
  morphisme de $S$-champs $[Y|G]\rightarrow \VECT S$, où $ \VECT S$
  est le champ des faisceaux localement libres (de rang fini) sur $S$.

  On note $G\Vect(Y)$ la catégorie des $G$-faisceaux localement libres sur $Y$.

\end{defi}

\begin{prop}
  \label{corrloc}
  Soit $X$ un schéma noethérien,
  $\pi : Y\rightarrow X$ un revêtement cyclique uniforme de degré
  $r$, $\mathcal F$ un $\boldsymbol{\mu}_r$-faisceau localement
  libre sur $Y$. Soit $D$ le diviseur de
  branchement et $E$ le diviseur de ramification de $\pi$.
  Alors l'association
\xymatrix@R=2pt{
  &&(\frac{1}{r}\mathbb Z)^{op}\ar[r] &\Vect(X)\\
   &&\frac{l}{r} \ar[r] & \pi_*^{\boldsymbol{\mu}_r}(\mathcal F\otimes_{\mathcal
    O_Y} \mathcal O_Y(-lE))}

  définit de manière naturelle un $r$-fibré parabolique sur $(X,D)$.
\end{prop}

\begin{rem}

  \begin{enumerate}

\item Ce type d'association est du à I.Biswas (\cite{Biswasorb}).

  \item
    Posons $\mathcal E_{\cdot}=\pi_*^{\boldsymbol{\mu}_r}(\mathcal
  F\otimes_{\mathcal O_Y} \mathcal O_Y(-\cdot r E))$. Pour ne pas
  alourdir la proposition nous n'avons défini $\mathcal E_{\cdot}$ que
  sur les objets. La définition sur les flèches est claire : pour
  $l\leq l'$, le morphisme $\mathcal E_{\frac{l'}{r}} \rightarrow \mathcal E_{\frac{l}{r}}$
  est induit par la multiplication par $t^{l'-l}$, où $t:\mathcal O_Y\rightarrow \mathcal O_Y(E)$ est la section canonique. Le morphisme
  $j:\mathcal E_{\cdot}\otimes_{\mathcal O_X}\mathcal O_X(-D)\simeq
  \mathcal E[1]_\cdot$ est donné par la relation $\pi^*D=rE$ et la
  formule de projection pour les $\boldsymbol{\mu}_r$-faisceaux. Enfin
  le fait que $j$ vérifie la condition de commutativité de la
  définition \ref{defYok} vient de la relation $\pi^* s=t^r$.

\item Nous verrons que cette opération définit un foncteur
  $\boldsymbol{\mu}_r\Vect(Y)\rightarrow \Par_{\frac{1}{r }}(X,D)$
  qui est en fait une équivalence de catégories tensorielles.

  \end{enumerate}
  \end{rem}

\begin{proof}[Démonstration de la proposition \ref{corrloc}]
  On commence par remarquer que
  pour $\mathcal F$ un ${\boldsymbol{\mu}_r}$-faisceau
 localement libre sur $Y$, $\pi_*^{\boldsymbol{\mu}_r}\mathcal F$
 est localement libre sur $X$ : en effet ${\boldsymbol{\mu}_r}$ étant diagonalisable c'est un facteur direct de
 $\pi_*\mathcal F$, qui est localement libre car $\mathcal F$ l'est,
  $\pi$ est fini et plat, et $X$ est noethérien.

Il reste à voir que  $l\leq l' < l+r $,
$\coker(\mathcal E_{\frac{l'}{r}} \rightarrow \mathcal E_{\frac{l}{r}})$
  est localement libre comme $\mathcal O_D$-module. Cela résulte
  clairement du lemme suivant :

  \begin{lem}
    Pour tout $\boldsymbol{\mu}_r$-faisceau  $\mathcal F$ localement
  libre sur $Y$, pour tout $0\leq l< r$, le
  faisceau  $\pi_*^{\boldsymbol{\mu}_r}(\mathcal
  F\otimes_{\mathcal O_Y} \mathcal O_{lE})$ est localement libre comme
  $\mathcal O_D$-module.
  \end{lem}
\begin{proof}
  La propriété à démontrer est locale sur $X$ et on peut donc supposer
  $X=\Spec R$, $Y=\Spec R[t]/(t^r - s)$ comme dans la définition \ref{RCU}.
On considère le diagramme commutatif suivant :
\begin{center}
\xymatrix@R=2pt@L=4pt{
  &&&lE \ar[ddd]_p\ar[r]^j & Y\ar[ddd]^{\pi}\\
  &&&&\\
  &&&&\\
   &&&D\ar[r]_i & X}
\end{center}
On a $\pi_*^{\boldsymbol{\mu}_r}(\mathcal
  F\otimes_{\mathcal O_Y} \mathcal O_{lE})=
  \pi_*^{\boldsymbol{\mu}_r}(j_*j^*\mathcal
  F)=i_*(p_*^{\boldsymbol{\mu}_r}j^*\mathcal F)$. Il s'agit donc de
  voir que $p_*^{\boldsymbol{\mu}_r}j^*\mathcal F$ est localement
  libre sur $D$. Or c'est un facteur direct de
  $p_*j^*\mathcal F$, qui est localement libre sur $D$, car
  $p: lE=\Spec \frac{R}{s}[t]/t^l \rightarrow \Spec \frac{R}{s}=D$ est
  fini et plat, $D$ est noethérien, et  $j^*\mathcal F$ est localement
  libre sur $lE$.
\end{proof}
\end{proof}

\section{Champ des racines $r$-ièmes d'un faisceau inversible muni d'une section}

On fixe un entier $r\geq 1$, un schéma noethérien
    $X\rightarrow\Spec \mathbb Z[r^{-1}]$, $\mathcal L$ un faisceau
    inversible sur $X$, $s\in H^0(X, \mathcal L)$.

\subsection{Définition}

\begin{defi}[\cite{Cadman}]

  On définit une catégorie $X_{(\mathcal L,s,r)}$ fibrée en groupoïdes
  sur $\Spec \mathbb Z[r^{-1}]$ par
\begin{enumerate}
\item  les objets de $X_{(\mathcal L,s,r)}$ sont les quadruplets
    $(f,\mathcal M, t, \phi)$ où $f:S\rightarrow X$ est un morphisme
    de schémas, $\mathcal M$ un faisceau inversible sur $S$,
    $t\in H^0(S, \mathcal M)$, et $\phi : \mathcal M^{\otimes
    r} \simeq f^*\mathcal L$ un isomorphisme tel que $\phi(t^{\otimes
    r})=f^*s$,
\item les morphismes de $(f,\mathcal M, t, \phi)$ (au dessus de
  $S$) à $(g,\mathcal N, u, \psi)$ (au dessus de
  $T$) sont les couples $(h, \rho)$, où $h:S\rightarrow T$ est un
  morphisme de schémas tel que $g\circ h =f$, et $\rho : \mathcal
  M\simeq h^*\mathcal N$ est un isomorphisme de faisceaux tel que
  $\rho(t)=h^*(u)$ et le diagramme suivant commute :
\begin{center}
\xymatrix@R=2pt@L=4pt{
  &&&\mathcal M^{\otimes r} \ar[ddd]_{\phi}\ar[r]^{\rho^{\otimes r}}
  & h^*\mathcal N^{\otimes r} \ar[ddd]^{h^*\psi}\\
  &&&&\\
  &&&&\\
   &&&f^*\mathcal L\ar@{=}[r] & h^*g^*\mathcal L }
\end{center}

\end{enumerate}
  \end{defi}

\begin{thm}[\cite{Cadman}, Theorem 2.2]
$X_{(\mathcal L,s,r)}$ est un champ de Deligne-Mumford
  (pour la topologie étale) sur $\Spec \mathbb Z[r^{-1}]$.
 \end{thm}

On appellera $X_{(\mathcal L,s,r)}$  le \emph{champ des racines $r$-ièmes} de
  $(\mathcal L,s)$.

\begin{prop}[\cite{Cadman}, Proposition 2.4]
Si $f:X'\rightarrow X$ est un morphisme de schémas alors
il existe un isomorphisme naturel
$X'\times_X X_{(\mathcal L,s,r)}\simeq X'_{(f^*\mathcal L,f^*s,r)}$.
\end{prop}

\begin{thm}[\cite{Cadman}, version 1, Proposition 3.2]\label{root}
Supposons qu'il existe un faisceau inversible $\mathcal N$ sur $X$
et un isomorphisme $\psi :\mathcal N^{\otimes r}\simeq \mathcal L$.

Si $Y=\SPEC (\frac{\Sym(\mathcal N^{\vee})}{{\mathcal N^{\vee}}^{\otimes r}\simeq \mathcal L^{\vee}})$, où la structure d'algèbre est définie par
$\psi^{\vee}$ et $s^{\vee}$, alors il existe un isomorphisme naturel
$[Y| {\boldsymbol{\mu}_r}]\simeq  X_{(\mathcal L,s,r)}$.
\end{thm}

\begin{proof}
  On rappelle brièvement le principe de la preuve.

  Soit $\pi : Y\rightarrow X$ le morphisme structurel, $u:\mathcal
  O_Y\rightarrow \pi^*\mathcal N$
  la section canonique. Alors $\xi:(\pi, \pi^*\mathcal N, u, \pi^*
  \psi)$ définit un objet de $X_{(\mathcal L,s,r)}(Y)$.

Donné un objet $(p,h)$ de $[Y| {\boldsymbol{\mu}_r}](S)$, composé
d'un ${\boldsymbol{\mu}_r}$-torseur $p:T\rightarrow S$, et d'un
morphisme ${\boldsymbol{\mu}_r}$-équivariant $h:T\rightarrow Y$, on
dispose donc d'un objet ${\boldsymbol{\mu}_r}$-équivariant $h^*\xi$.
Comme $X_{(\mathcal L,s,r)}$ est un champ celui-ci se descend en un
objet $\eta$ de $X_{(\mathcal L,s,r)}(S)$.

Réciproquement fixons un objet $(f,\mathcal M, t, \phi)$ de
$X_{(\mathcal L,s,r)}(S)$. L'isomorphisme $\mathcal M^{\otimes
  r}\simeq f^*\mathcal L\simeq f^*\mathcal N^{\otimes
  r}$ déduit de $\phi$ et $f^*\psi$ définit un
  ${\boldsymbol{\mu}_r}$-torseur sur $S$

$$T=\SPEC (\frac{\Sym(\mathcal M \otimes f^*\mathcal N^{\vee})}{\mathcal M^{\otimes r}\simeq f^*\mathcal N^{\otimes r}})$$

  et $t:\mathcal O_S\rightarrow \mathcal M$ induit un morphisme
  ${\boldsymbol{\mu}_r}$-équivariant $T\rightarrow Y$. On dispose donc
  d'un objet de $[Y| {\boldsymbol{\mu}_r}](S)$, et on vérifie que ces
  deux constructions définissent des équivalences de catégories
  réciproques l'une de l'autre.

\end{proof}

\begin{cor}
  \label{strucloc}

  Pour tout point $x$ de $X$, il existe un voisinage affine $U =  \Spec R$
de  $x$ dans $X$, un élément $\sigma \in R$ et un isomorphisme naturel de $U$-champs :
$$[\Spec (\frac{R[\tau]}{\tau^r - \sigma})| {\boldsymbol{\mu}_r}]\simeq U_{(\mathcal
  L_{|U},s_{|U},r)}$$
\end{cor}

\begin{cor}
\label{pipropre}
Le morphisme canonique $\pi : X_{(\mathcal L,s,r)}\rightarrow X$ est fini. 
\end{cor}

\begin{proof}
On doit montrer que $\pi$ est quasi-fini et propre.

Commençons par la quasi-finitude :  d'après \cite{Cadman}, on dispose d'un
atlas étale $U\rightarrow X_{(\mathcal L,s,r)}$ de type fini sur $X$, donc $\pi$ est de type fini. Donné un point géométrique $\Spec \Omega\rightarrow X$, 
on sait que $\Spec \Omega\times_ X X_{(\mathcal L,s,r)}\simeq 
\Spec \Omega_{(\mathcal L_{|\Spec \Omega},s_{|\Spec \Omega},r)}$, et d'après le théorème \ref{root}, ce champ admet bien un atlas étale fini sur $\Spec \Omega$. Donc 
$\pi$ est quasi-fini.

Prouvons que $\pi$ est propre. On s'assure aisément que $\pi$ vérifie le critère valuatif usuel (\cite{VistInv}) pour les morphismes universellement fermés.
Reste à voir que $\pi$ est séparé, ou de manière équivalente que le morphisme diagonal $X_{(\mathcal L,s,r)}\rightarrow X_{(\mathcal L,s,r)}\times_X X_{(\mathcal L,s,r)}$ est propre. Soit $(X_i\rightarrow X)$ un revêtement Zariski de $X$ trivialisant $\mathcal L$, et $\mathcal L_i$, $s_i$ les restrictions. D'après (\cite{VistDesc}), Proposition 2.36, il suffit de montrer que  
${X_i}_{(\mathcal L_i,s_i,r)}\rightarrow {X_i}_{(\mathcal L_i,s_i,r)}\times_{X_i}{X_i}_{(\mathcal L_i,s_i,r)}$ est propre. On peut donc supposer $\mathcal L$ trivial, et 
d'après le théorème \ref{root}, que 
$X_{(\mathcal L,s,r)}\simeq [Y| {\boldsymbol{\mu}_r}]$, où $Y\rightarrow X$
est fini. 

Le changement de base du morphisme diagonal 
$[Y| {\boldsymbol{\mu}_r}]\rightarrow
[Y| {\boldsymbol{\mu}_r}]\times_X[Y| {\boldsymbol{\mu}_r}]$ par l'atlas étale
$Y \times_X Y$ de $[Y| {\boldsymbol{\mu}_r}]\times_X[Y| {\boldsymbol{\mu}_r}]$
donne le morphisme $(a,pr_2): {\boldsymbol{\mu}_r}\times Y\rightarrow
Y\times_X Y$, où $a$ désigne le morphisme d'action.
Or celui-ci est fini (et donc propre), car fini composé avec 
$pr_2:  Y\times_X Y\rightarrow Y$, qui est séparé car obtenu par changement de
base à partir du morphisme affine $Y\rightarrow X$. 
\end{proof}
 
\begin{cor}
$X_{(\mathcal L,s,r)}$ admet $X$ comme espace des modules grossier.
\end{cor}

\begin{proof}
C'est un conséquence immédiate du fait que $X_{(\mathcal L,s,r)}\rightarrow X$ est séparé, donc admet un espace des modules grossier (voir \cite{KM}), et que $X$ est un quotient catégorique pour $X_{(\mathcal L,s,r)}$ (i.e. est universel pour les morphismes de $X_{(\mathcal L,s,r)}$ vers un schéma) comme montré dans \cite{Cadman}, \S 2.4.
\end{proof}

  \begin{cor}
 \label{quotunif}
  Soit $\pi : Y\rightarrow X$ un revêtement cyclique uniforme de degré
  $r$. Soit $D$ le diviseur (resp. $E$) de branchement (resp. de
   ramification), et $s$ (resp. $t$) la section canonique de
   $\mathcal L =\mathcal O_X(D)$ (resp.  $\mathcal N=\mathcal O_Y(E)$).
Soit de plus $\phi : \mathcal N^{\otimes r} \simeq \pi^* \mathcal L$
   le morphisme découlant de l'égalité $\pi^* D =rE$.

 Alors le morphisme $Y \rightarrow  X_{(\mathcal L,s,r)}$ défini par
 $(\pi,\mathcal N, t,\phi)$ induit un isomorphisme de $X$-champs
 $[Y| {\boldsymbol{\mu}_r}]\simeq  X_{(\mathcal L,s,r)}$
\end{cor}

  \begin{proof}

Cela découle du théorème \ref{root} et de la remarque \ref{remRCU}.
\end{proof}

 \begin{cor}
\label{pilisse}
Supposons que $s:\mathcal O_X\rightarrow \mathcal L$ est un monomorphisme, et notons $D=\Div(s)$. Supposons de plus donné un morphisme lisse $X\rightarrow S_0$, où $S_0\rightarrow \Spec \mathbb Z[r^{-1}]$ est une base, tel que $D$ est lisse sur $S_0$. Alors $X_{(\mathcal L,s,r)}$ est lisse sur $S_0$.
\end{cor}

\begin{proof}
La lissité est une propriété locale, il s'agit donc de voir que 
$X_{(\mathcal L,s,r)}$ admet un atlas lisse sur $S_0$. Quitte à recouvrir $X$ par des ouverts trivialisant $\mathcal L$, on peut se ramener, d'après le théorème \ref{root}, au cas où $X_{(\mathcal L,s,r)}\simeq [Y| {\boldsymbol{\mu}_r}]$, où $Y\rightarrow X$  est un revêtement cyclique uniforme. Mais alors 
$[Y| {\boldsymbol{\mu}_r}]$ admet $Y\rightarrow [Y| {\boldsymbol{\mu}_r}]$ 
comme atlas étale, et $Y\rightarrow S_0$ est lisse d'après \cite{AV}, Proposition 2.5. 

\end{proof}

\begin{lem}
\label{CartDisj}
  Soient $D,D'$ deux diviseurs de Cartier effectifs à supports disjoints,
 $\mathcal L =\mathcal O_X(D)$, $\mathcal L' =\mathcal O_X(D')$, et $s$, $s'$ 
les sections canoniques. Il existe un isomorphisme canonique
$X_{(\mathcal L,s,r)}\times_X X_{(\mathcal L',s',r)}\simeq 
X_{(\mathcal L\otimes_{\mathcal O_X} \mathcal L',s\otimes s',r)}$.
\end{lem}

\begin{proof}
  Le morphisme canonique $X_{(\mathcal L,s,r)}\times_X X_{(\mathcal
    L',s',r)}\rightarrow X_{(\mathcal L\otimes_{\mathcal O_X} \mathcal
    L',s\otimes s',r)}$ est un isomorphisme au dessus de $U=X-\supp(D)$ et
$U'=X-\supp(D')$, qui par hypothèse recouvrent $X$. 

\end{proof}
  \subsection{Faisceaux inversibles sur ce champ}

  Sur $X_{(\mathcal L,s,r)}$ on dispose d'un faisceau inversible canonique
$\mathcal N$  défini par $(f,\mathcal M, t, \phi)\rightarrow \mathcal
M$. Il est muni d'une section canonique $u : \mathcal
O_{X_{(\mathcal L,s,r)}}\rightarrow \mathcal N$ définie localement par
$t:\mathcal O_S \rightarrow \mathcal M$.

De plus, il existe un isomorphisme canonique $\psi : \mathcal N^{\otimes
  r}\simeq \pi^* \mathcal L$, défini localement par $\phi: \mathcal M^{\otimes
  r}\simeq f^* \mathcal L$, il vérifie $\psi(u^ {\otimes
  r})=\pi^* s$.

  \begin{lem}
  \label{PartEnt}
  On suppose que $\Div s$ est un diviseur de Cartier effectif.
  Pour tout entier $l$ le morphisme naturel :
  $$\mathcal L^{\otimes[l/r]}\rightarrow \pi_*(\mathcal N^{\otimes l}) $$
est un isomorphisme.
\end{lem}

\begin{proof}
  D'après la formule de projection, on a pour tout entiers $a,b$ :
  $\pi_*(\mathcal N^{\otimes ar+b})\simeq \mathcal L^{\otimes a}
  \pi_*(\mathcal N^{\otimes b})$, on
  peut donc se contenter de montrer la formule pour $r$ entiers
  consécutifs.

  Soit $0< l \leq r$. La suite exacte

$$0\rightarrow \mathcal N^{\otimes-l}\rightarrow\mathcal O_{X_{(\mathcal L,s,r)}}
\rightarrow  \mathcal O_{X_{(\mathcal L,s,r)}}/\mathcal
N^{\otimes-l}\rightarrow 0$$

correspond localement à une suite de ${\boldsymbol{\mu}_r}$-faisceaux

$$0\rightarrow\frac{t^lR[t]}{t^r-s}\rightarrow
\frac{R[t]}{t^r-s}\rightarrow \frac{R/s[t]}{t^l}\rightarrow 0$$

En prenant les fixes il vient
$$0\rightarrow(\frac{t^lR[t]}{t^r-s})^{\boldsymbol{\mu}_r}\rightarrow
R\rightarrow R/s\rightarrow 0$$

Comme $s$ est l'équation locale définissant $\mathcal L$, on en déduit
que $\pi_* (\mathcal N^{\otimes-l})=\mathcal L^{\otimes-1}$, d'où la conclusion.

\end{proof}

\subsection{Faisceaux localement libres sur ce champ}

Pour $\mathfrak X$ un champ de Deligne-Mumford sur le schéma $X$,
un faisceau localement libre sur $\mathfrak X$ est par définition un
morphisme de champs $\mathfrak X\rightarrow \VECT X$, où $\VECT X$ est
le champ des faisceaux localement libres sur $X$.
On note $\Vect \mathfrak X$ la catégorie des faisceaux localement
libres sur $\mathfrak X$.

\subsubsection{Situation locale}

 \begin{prop}
  \label{StrucLoc}
  On suppose que $\Div s$ est un diviseur de Cartier effectif.
Soit $\mathcal F$ un faisceau localement libre sur $X_{(\mathcal L,s,r)}$.
Pour tout point $x$ de $X$, il existe un voisinage $U$ (pour la topologie de Zariski) de $x$ tel que
$\mathcal F_{|\pi^{-1}U}$ est une somme de faisceaux inversibles.
\end{prop}

 \begin{proof}
 D'après le corollaire \ref{strucloc}, on peut identifier
 $X_{(\mathcal L,s,r)}\rightarrow X$ à $[(\Spec R[t]/(t^r - s))|
 {\boldsymbol{\mu}_r}]\rightarrow \Spec R$, pour $R$ un anneau, $s$ un
 élément de $R$ ne divisant pas zéro. Soit $R'= R[t]/(t^r - s)$.
La donnée de $\mathcal F$ équivaut à celle d'un $R'$-module libre $M$,
 avec une $\mathbb Z/r$-graduation compatible avec celle de $R'$ et
 l'action de $R'$ sur $M$, on la notera $M=M_0\oplus M_1\oplus \cdots
 \oplus M_{r-1}$.

On peut de plus supposer que $R$
 est un anneau local, d'idéal maximal $\mathfrak m$.
Si $s\notin \mathfrak m$ il n'y a rien à démontrer, puisque alors
 ${\boldsymbol{\mu}_r}$ agit librement sur $\Spec R[t]/(t^r - s)$ et
 donc $[(\Spec R[t]/(t^r - s))|{\boldsymbol{\mu}_r}]\simeq \Spec R$.

 On peut donc supposer que $s\in \mathfrak m$. Alors $R'= R[t]/(t^r -
 s)$ est aussi local. En effet, si $\mathfrak n$ est un idéal maximal
 de $R'$, alors $R'/R$ étant finie, on a $R\cap \mathfrak n=\mathfrak
 m$, et donc $s\in \mathfrak n$. Donc $t\in \mathfrak n$,
 et $\mathfrak n$ peut-être vu comme
 un idéal maximal de $R'/t$, comme $R'/t\simeq R/s$ est local, $ \mathfrak n$
 est unique.

 Notons que pour $l$ entier, on a un isomorphisme naturel
 de $R'$-modules gradués $R'[l]\simeq t^lR'$, en particulier
 pour $l=1$ on obtient un isomorphisme :

 $$\frac{M}{tM}\simeq \oplus_j\frac{M_j}{M_{j+1}}$$

 de $R'/t$-modules libres, respectant la graduation. Pour chaque $j$,
 choisissons un $n_j$-uplet $(e_1^j,\cdots,e_{n_j}^j)$ dans $M_j$ tel que
 son image forme une base de $\frac{M_j}{M_{j+1}}$, et soit
 $(e_1,\cdots,e_n)$ le $n$-uplet obtenu par concaténation.

 Comme  $M_j\simeq \Hom_{R'}( R'[j], M)$ (morphismes de $R'$-modules
 gradués) ce $n$-uplet définit un morphisme de $R'$-modules gradués

 $$\oplus_j R'[j]^{\oplus n_j}\rightarrow M$$

 Comme $R'$ est local, et $t\in \mathfrak n$, le lemme de Nakayama
 implique que ce morphisme est surjectif. Comme les deux $R'$-modules
 sont libres de même rang ce doit être un isomorphisme de
 $R'$-modules, et donc un isomorphisme de $R'$-modules gradués.

 \end{proof}

\subsubsection{Situation globale}

\begin{thm}
\label{thmprinc}
   Soit $r$ un entier naturel,
   $X\rightarrow\Spec \mathbb Z[r^{-1}]$ un schéma noethérien,
 $D$ un diviseur de Cartier effectif, $\mathcal L=\mathcal O_X(D)$, $s$ la
 section canonique, alors le foncteur $F$ :

\xymatrix@R=2pt{
  &&\Vect (X_{(\mathcal L,s,r)})\ar[r]&\Par_{\frac{1}{r}} (X,D)& \\
  &&\mathcal F\ar[r] & ((\frac{1}{r}\mathbb Z)^{op}  \ar[r]
  &\Vect(X))\\
&&& \frac{l}{r}\ar[r] & \pi_* (\mathcal N^{\otimes -l}\otimes_{\mathcal O_{X_{(\mathcal L,s,r)}}} \mathcal F)}

est une équivalence de catégories tensorielles.
\end{thm}

\subsection{Preuve}

\subsubsection{$F$ est bien défini}

\begin{lem}
 Pour tout faisceau $\mathcal F$ localement libre sur $X_{(\mathcal
  L,s,r)}$, le foncteur $ \frac{l}{r}\rightarrow  \pi_* (\mathcal
  N^{\otimes -l}\otimes_{\mathcal O_{X_{(\mathcal L,s,r)}}} \mathcal F)$
  définit un $r$-fibré parabolique sur $(X,D)$.
\end{lem}

\begin{proof}

Posons $\mathcal E_{\frac{l}{r}}= \pi_* (\mathcal
  N^{\otimes-l}\otimes_{\mathcal O_{X_{(\mathcal L,s,r)}}} \mathcal F)$.
Il s'agit de voir que $\mathcal E_{\frac{l}{r}}$ est localement libre
  sur $X$
et que pour $l\leq l' < l+r $, $\coker(\mathcal E_{\frac{l'}{r}} \rightarrow
\mathcal E_{\frac{l}{r}})$ est localement libre comme $\mathcal O_D$-module.
Ce sont des affirmations locales sur $X$. En trivialisant $\mathcal L$
on se ramène au cas où $X_{(\mathcal L,s,r)}$ est un champ quotient
issu d'un revêtement cyclique uniforme, comme dans le corollaire
\ref{strucloc}, et on peut alors appliquer la proposition \ref{corrloc}.

\end{proof}

\subsubsection{Définition d'un foncteur quasi-inverse $G$}

Pour les bouts, voir l'appendice \ref{Bout}.

\begin{prop}
  \label{BoutExist}
  Soit $\mathcal E_\cdot \in \obj \Par_{\frac{1}{r}} (X,D)$.

  Le bout universel $\int^{\frac{1}{r}\mathbb Z}\mathcal N^{\otimes l}
  \otimes_{\mathcal O_{X_{(\mathcal L,s,r)}}}\pi^*(\mathcal
  E_{\frac{l}{r}})$ existe dans $\Vect (X_{(\mathcal L,s,r)})$.
\end{prop}

\begin{proof}

On commence par remarquer que la question a un sens : en effet
$(\frac{l}{r},\frac{l'}{r}) \rightarrow \mathcal N^{\otimes l'}
  \otimes_{\mathcal O_{X_{(\mathcal L,s,r)}}}\pi^*(\mathcal
  E_{\frac{l}{r}})$ définit un foncteur de variance mixte
$(\frac{1}{r}\mathbb Z)^{op}\times \frac{1}{r}\mathbb Z \rightarrow
  \Vect (X_{(\mathcal L,s,r)})$, la contravariance en la première variable
  venant de la structure de $r$-faisceau parabolique de $\mathcal
  E_\cdot$, et la covariance en la seconde venant de la multiplication
  par une puissance convenable de $u : \mathcal
O_{X_{(\mathcal L,s,r)}}\rightarrow \mathcal N$.

Un bout universel est une colimite (\cite{MacLane}), et on peut donc
la construire localement : en effet, le caractère universel
des limites locales rend les conditions de recollement
automatiques. Plus précisément, on utilise le lemme folklorique suivant, qu'on
donne ici sans démonstration.

\begin{lem}
Soit $\mathfrak X$ un champ algébrique, $\mathcal E :I \rightarrow \Vect \mathfrak X$ un diagramme, $(\mathfrak X_\alpha \rightarrow \mathfrak X)_\alpha $ un
recouvrement par des ouverts.
Si $$\lim_{\stackrel{\rightarrow}{I}} {\mathcal E_i}_{|\mathfrak X_\alpha}$$ existe dans
$\Vect \mathfrak X_\alpha$ pour tout $\alpha$, alors
 $$\lim_{\stackrel{\rightarrow}{I}} {\mathcal E_i}$$
  existe dans $\Vect \mathfrak X$.
\end{lem}

Grâce au théorème \ref{LocStrPar}, en raisonnant localement sur $X$,
on peut donc se contenter de montrer
l'existence de $\int^{\frac{1}{r}\mathbb Z}\mathcal N^{\otimes l}
  \otimes_{\mathcal O_{X_{(\mathcal L,s,r)}}}\pi^*(\mathcal
  E_{\frac{l}{r}})$ lorsque $\mathcal E_\cdot =
  \underline{\mathcal O_X}_\cdot [\frac{m}{r}]$, i.e. $\mathcal
  E_\cdot$ est un décalé du faisceau structurel muni de la structure
  spéciale. Or $\underline{\mathcal O_X}_\cdot [\frac{m}{r}]$ est
  défini sur les objets par
  $(\underline{\mathcal O_X}_\cdot
  [\frac{m}{r}])_{\frac{l}{r}}=\mathcal L ^{\otimes [-\frac{l+m}{r}]}$ et sur
  les flèches par multiplication par une puissance convenable de $s
  :\mathcal O_X \rightarrow  \mathcal L$. Il
  reste à prouver le lemme suivant :

  \begin{lem}
    \label{reduc}
$$\int^{\frac{1}{r}\mathbb Z}\mathcal N^{\otimes l}
  \otimes_{\mathcal O_{X_{(\mathcal L,s,r)}}}\pi^*\mathcal L^{\otimes [-\frac{l+m}{r}]}
\simeq \mathcal N^{\otimes -m}$$
  \end{lem}

\begin{proof}
  L'isomorphisme $\psi : \mathcal N^{\otimes
  r}\simeq \pi^* \mathcal L$, et la relation $\psi(u^ {\otimes
  r})=\pi^* s$, montrent que les foncteurs $(l,l')\rightarrow \mathcal N^{\otimes l'}
  \otimes_{\mathcal O_{X_{(\mathcal L,s,r)}}}\pi^*\mathcal L
  ^{\otimes [-\frac{l+m}{r}]}$ et $(l,l')\rightarrow \mathcal N^{\otimes l'+
  r[-\frac{l+m}{r}]}$ sont isomorphes (les flèches du second sont
  définies par multiplication par une puissance convenable de $u : \mathcal
O_{X_{(\mathcal L,s,r)}}\rightarrow \mathcal N$).
Or pour tout $l$ entier on a : $ -m-r<l+r[-\frac{l+m}{r}]\leq -m$, et la
  multiplication par $u$ fournit un morphisme $\mathcal N^{\otimes l+
  r[-\frac{l+m}{r}]}\rightarrow \mathcal N^{\otimes -m}$, évidemment dinaturel. Comme
  $l+r[-\frac{l+m}{r}]$ prend la valeur $-m$ lorsque $l$ varie, ce
  bout est universel, d'où
$\int^{\frac{1}{r}\mathbb Z}\mathcal N^{\otimes l+
  r[-\frac{l+m}{r}]}\simeq \mathcal N^{\otimes -m}$, et la conclusion.
\end{proof}

\end{proof}

La proposition \ref{BoutExist} permet de poser
$G(\mathcal E_\cdot)=\int^{\frac{1}{r}\mathbb Z}\mathcal N^{\otimes l}
  \otimes_{\mathcal O_{X_{(\mathcal L,s,r)}}}\pi^*(\mathcal
  E_{\frac{l}{r}})$. La dinaturalité du bout rend cette expression
  fonctorielle en $\mathcal E_\cdot$, si bien qu'on a défini un
  foncteur $G :\Par_{\frac{1}{r}} (X,D)  \rightarrow \Vect
  (X_{(\mathcal L,s,r)})$.

\subsubsection{Preuve de l'équivalence}

\boxed{G\circ F \simeq 1}

On doit construire un
isomorphisme $$\int^{\frac{1}{r}\mathbb Z}\mathcal N^{\otimes l}
  \otimes_{\mathcal O_{X_{(\mathcal L,s,r)}}}\pi^*
\pi_* (\mathcal N^{\otimes -l}\otimes_{\mathcal O_{X_{(\mathcal L,s,r)}}}
\mathcal F)\rightarrow \mathcal F$$
naturel en $\mathcal F\in \obj \Vect (X_{(\mathcal L,s,r)})$.

L'existence d'un morphisme naturel est donné par le lemme suivant.

\begin{lem}
Soit $adj : \Hom(\pi^* \mathcal E, \mathcal F)\simeq \Hom(\mathcal E,
\pi_*\mathcal F)$ le morphisme
d'adjonction. La famille de morphismes
$$1\otimes adj^{-1}(1) : \mathcal N^{\otimes l}
  \otimes_{\mathcal O_{X_{(\mathcal L,s,r)}}}\pi^*
\pi_* (\mathcal N^{\otimes -l}\otimes_{\mathcal O_{X_{(\mathcal L,s,r)}}}
\mathcal F) \rightarrow \mathcal F$$ est dinaturelle en $l\in \obj
\mathbb Z$.
\end{lem}

\begin{proof}
Il s'agit de voir la commutativité du diagramme suivant :

\begin{center}
\xymatrix@R=4pt{
&&\mathcal N^{\otimes l}
  \otimes\pi^*
\pi_* (\mathcal N^{\otimes -l}\otimes
\mathcal F)\ar[drr]_{1\otimes adj^{-1}(1)}  && \\
\mathcal N^{\otimes l'}
  \otimes\pi^*
\pi_* (\mathcal N^{\otimes -l}\otimes
\mathcal F)\ar[urr]^{u^{l-l'}\otimes 1}\ar[drr]^{1\otimes \pi^*\pi_*(u^{l-l'}\otimes 1)}&&& &\mathcal F\\
&&\mathcal N^{\otimes l'}
  \otimes\pi^*
\pi_* (\mathcal N^{\otimes -l'}\otimes
\mathcal F)\ar[urr]^{1\otimes adj^{-1}(1)}  &&}
\end{center}

ou de manière équivalente celle du diagramme :

\begin{center}
\xymatrix@R=2pt{
&&\mathcal N^{\otimes l-l'}
  \otimes\pi^*
\pi_* (\mathcal N^{\otimes -l}\otimes
\mathcal F)\ar[drr]_{1\otimes adj^{-1}(1)}  && \\
\pi^*\pi_* (\mathcal N^{\otimes -l}\otimes
\mathcal F)\ar[urr]^{u^{l-l'}}\ar[drr]^{\pi^*\pi_*(u^{l-l'}\otimes 1)}&&& &\mathcal N^{\otimes -l'}
  \otimes\mathcal F\\
&&\pi^*
\pi_* (\mathcal N^{\otimes -l'}\otimes
\mathcal F)\ar[urr]_{adj^{-1}(1)}  &&}
\end{center}

ce qui équivaut encore à celle du diagramme :

\begin{center}
\xymatrix@R=2pt{
&&\mathcal N^{\otimes -l}
  \otimes\mathcal F \ar[drr]^{u^{l-l'}\otimes 1}  && \\
\pi^*\pi_* (\mathcal N^{\otimes -l}\otimes
\mathcal F)\ar[urr]^{adj^{-1}(1)}\ar[drr]^{\pi^*\pi_*(u^{l-l'}\otimes 1)}&&& &\mathcal N^{\otimes -l'}
  \otimes\mathcal F\\
&&\pi^*
\pi_* (\mathcal N^{\otimes -l'}\otimes
\mathcal F)\ar[urr]_{adj^{-1}(1)}  &&}
\end{center}

Or $adj$ envoie ces deux flèches sur $\pi_*(u^{l-l'}\otimes 1)$,
d'où la conclusion.
\end{proof}

Vérifier que le morphisme naturel :
$\int^{\frac{1}{r}\mathbb Z}\mathcal N^{\otimes l}
  \otimes_{\mathcal O_{X_{(\mathcal L,s,r)}}}\pi^*
\pi_* (\mathcal N^{\otimes -l}\otimes_{\mathcal O_{X_{(\mathcal L,s,r)}}}
\mathcal F)\rightarrow \mathcal F$
est un isomorphisme est une question locale, et on peut donc supposer
que $X=\Spec R$, où $R$ est un anneau local. On peut aussi supposer 
que $\mathcal F$ est un faisceau inversible d'après la proposition \ref{StrucLoc}. 
Or $R$ étant local, $\Pic(X_{(\mathcal
  L,s,r)})$ est cyclique d'ordre $r$, et engendré par $\mathcal N$. On
peut donc prendre $\mathcal F=\mathcal N^{\otimes -m}$.

Finalement d'après le lemme \ref{PartEnt}
$\int^{\frac{1}{r}\mathbb Z}\mathcal N^{\otimes l}
  \otimes \pi^*\pi_* \mathcal N^{\otimes -l-m}\simeq
  \int^{\frac{1}{r}\mathbb Z}\mathcal N^{\otimes l} \otimes\pi^*\mathcal
  L^{\otimes [-\frac{l+m}{r}]}$ et on peut conclure à l'aide du lemme \ref{reduc}.

  \ \\
\boxed{F\circ G \simeq 1}

On doit montrer l'existence d'un inverse pour le morphisme naturel :

$$\pi_*(\mathcal N^{\otimes -m}\otimes_{\mathcal O_{X_{(\mathcal L,s,r)}}}
\int^{\frac{1}{r}\mathbb Z}\mathcal N^{\otimes l}
  \otimes_{\mathcal O_{X_{(\mathcal L,s,r)}}}\pi^*(\mathcal
E_{\frac{l}{r}}))\leftarrow \mathcal
E_{\frac{m}{r}}$$

Le foncteur $\mathcal N^{\otimes -l}\otimes_{\mathcal O_{X_{(\mathcal
      L,s,r)}}}\cdot $ étant une équivalence on a :

$$\pi_*(\mathcal N^{\otimes -m}\otimes_{\mathcal O_{X_{(\mathcal L,s,r)}}}
\int^{\frac{1}{r}\mathbb Z}\mathcal N^{\otimes l}
  \otimes_{\mathcal O_{X_{(\mathcal L,s,r)}}}\pi^*(\mathcal
E_{\frac{l}{r}}))\simeq
\pi_*(\int^{\frac{1}{r}\mathbb Z}\mathcal N^{\otimes l-m}
  \otimes_{\mathcal O_{X_{(\mathcal L,s,r)}}}\pi^*(\mathcal
E_{\frac{l}{r}}))$$

Comme ${\boldsymbol{\mu}_r}$ est diagonalisable,
$\pi_*$ est exact à droite, et commute aux sommes
finies, donc aux colimites finies, et on a donc :

$$\pi_*(\int^{\frac{1}{r}\mathbb Z}\mathcal N^{\otimes l-m}
  \otimes_{\mathcal O_{X_{(\mathcal L,s,r)}}}\pi^*(\mathcal
E_{\frac{l}{r}}))\simeq
\int^{\frac{1}{r}\mathbb Z}\pi_*(\mathcal N^{\otimes l-m}
  \otimes_{\mathcal O_{X_{(\mathcal L,s,r)}}}\pi^*(\mathcal
E_{\frac{l}{r}}))$$

Comme $\mathcal
E_{\frac{l}{r}}$ est localement libre, on a d'après la formule de
projection :

$$
 \int^{\frac{1}{r}\mathbb Z}\pi_*(\mathcal N^{\otimes l-m}
  \otimes_{\mathcal O_{X_{(\mathcal L,s,r)}}}\pi^*(\mathcal
E_{\frac{l}{r}}))\simeq \int^{\frac{1}{r}\mathbb Z}\pi_*(\mathcal N^{\otimes l-m})
  \otimes_{\mathcal O_X}\mathcal
E_{\frac{l}{r}}$$

Le lemme \ref{PartEnt} donne

$$
\int^{\frac{1}{r}\mathbb Z}\pi_*(\mathcal N^{\otimes l-m})
  \otimes_{\mathcal O_X}\mathcal
E_{\frac{l}{r}}\simeq \int^{\frac{1}{r}\mathbb Z}
\mathcal L^{\otimes [\frac{l-m}{r}]} \otimes_{\mathcal O_X}\mathcal
E_{\frac{l}{r}}$$

Finalement, $\mathcal E_\cdot$ est un fibré parabolique, donc muni
d'un isomorphisme $\mathcal L^{\otimes-1}\otimes _{\mathcal O_X}\mathcal
E_\cdot\simeq \mathcal E_\cdot[1]$, qui induit

$$\int^{\frac{1}{r}\mathbb Z}
\mathcal L^{\otimes [\frac{l-m}{r}]} \otimes_{\mathcal O_X}\mathcal
E_{\frac{l}{r}}\simeq \int^{\frac{1}{r}\mathbb Z} \mathcal
E_{\frac{l}{r}-[\frac{l-m}{r}]}$$

Finalement le fait que le morphisme naturel

$$ \int^{\frac{1}{r}\mathbb Z} \mathcal
E_{\frac{l}{r}-[\frac{l-m}{r}]}\leftarrow\mathcal E_{\frac{m}{r}}$$

est un isomorphisme résulte du fait que
  ${\frac{l}{r}-[\frac{l-m}{r}]}$ prend la valeur ${\frac{m}{r}}$
  lorsque $l$ varie, et de l'unicité des bouts universels.

\subsubsection{Preuve du caractère tensoriel}

D'après \cite{Saav}, I 4.4 on peut se contenter de montrer que le
foncteur $G$ est compatible avec le produit tensoriel. On commence par 
réinterpréter le produit tensoriel des fibrés paraboliques en termes de bout. 

\begin{lem}
Soient $\mathcal E_\cdot$, $\mathcal E'_\cdot$ deux objets de
$\Par_{\frac{1}{r}} (X,D)$. Alors

$$(\mathcal E_\cdot\otimes\mathcal E'_\cdot)_{\frac{l}{r}}
\simeq\int^{\frac{1}{r}\mathbb Z} \mathcal
E_{\frac{m}{r}}\otimes_{\mathcal O_X} \mathcal
E'_{\frac{l-m}{r}}$$

\end{lem}

\begin{proof}
  D'après la remarque \ref{remfibrpar}, on peut plonger 
$\Par_{\frac{1}{r}} (X,D)$ comme une sous-catégorie pleine de la catégorie de 
foncteurs $\Func((\frac{1}{r}\mathbb Z)^{op},\MOD X)$, catégorie dans laquelle on peut donc calculer le produit tensoriel $\mathcal E_\cdot\otimes\mathcal E'_\cdot$.
La formule de convolution donnée résulte alors de \cite{Day}, \S 4.
\end{proof}
Il s'agit de comparer

$$G(\mathcal E_\cdot)\otimes G(\mathcal E'_\cdot)=
\int^{\frac{1}{r}\mathbb Z}\mathcal N^{\otimes l}
  \otimes\pi^*(\mathcal
  E_{\frac{l}{r}})\otimes
  \int^{\frac{1}{r}\mathbb Z}\mathcal N^{\otimes l'}
  \otimes\pi^*(\mathcal
  E'_{\frac{l'}{r}})
  \simeq
  \int^{\frac{1}{r}\mathbb Z}\int^{\frac{1}{r}\mathbb Z}\mathcal N^{\otimes l+l'}\otimes \pi^*(\mathcal
  E_{\frac{l}{r}}\otimes \mathcal
  E'_{\frac{l'}{r}})$$

  et
 $$G(\mathcal E_\cdot \otimes \mathcal E'_\cdot)=
  \int^{\frac{1}{r}\mathbb Z}\mathcal N^{\otimes l}
  \otimes\pi^*((\mathcal E_\cdot\otimes\mathcal
  E'_\cdot)_{\frac{l}{r}})\simeq
\int^{\frac{1}{r}\mathbb Z}\mathcal N^{\otimes l}
  \otimes\pi^*( \int^{\frac{1}{r}\mathbb Z} \mathcal
E_{\frac{m}{r}}\otimes_{\mathcal O_X} \mathcal
E'_{\frac{l-m}{r}})
  $$

$\pi^*$ est adjoint à gauche, donc commute aux colimites, et donc

$$G(\mathcal E_\cdot \otimes \mathcal E'_\cdot)\simeq \int^{\frac{1}{r}\mathbb Z}\int^{\frac{1}{r}\mathbb Z} \mathcal N^{\otimes l}
  \otimes\pi^*( \mathcal
E_{\frac{m}{r}}\otimes_{\mathcal O_X} \mathcal
E'_{\frac{l-m}{r}})
  $$

On peut donc conclure en appliquant le théorème de Fubini pour les
bouts universels (théorème \ref{Fubini} de l'appendice \ref{Bout}).

\section{Degrés}

Dans toute cette partie, on fixe $r\geq 1$ un entier naturel,
$X$ une variété projective lisse sur un corps $k$ algébriquement clos de
caractéristique première à $r$. On note $\mathcal O(1)$ un faisceau inversible très ample sur $X$. On se fixe de plus $D$ un diviseur de Cartier effectif sur $X$, lisse sur $k$. On note $n$ la dimension de $X$.

\subsection{Définitions}

\subsubsection{Degré parabolique}

\begin{defi}[\cite{MY}]\label{degpar}
Soit $\mathcal E_\cdot$ un objet de $\Par_{\frac{1}{r}} (X,D)$, de rang $\rho$.
On définit sa \emph{caractéristique d'Euler parabolique} par :
$$\chi_{par}(\mathcal E_\cdot)=\frac{1}{r}\sum_{l=1}^{r}\chi(X,\mathcal E_\frac{l}{r})$$
Son \emph{degré parabolique} $\deg_{par}(\mathcal E_\cdot)$ est par définition le nombre rationnel défini par:
$$\deg_{par}(\mathcal E_\cdot)=(n-1)!\times \{\text{coefficient de } m^{n-1} \text{dans}\;  \chi_{par}(\mathcal E_\cdot(m))- \chi_{par}(\underline{\mathcal O_X}_\cdot^{\oplus \rho}(m)) \}$$
\end{defi}

\subsubsection{Degré sur un champ de Deligne-Mumford projectif}

\begin{defi}
Soit $\mathfrak X\rightarrow S$ un champ de Deligne-Mumford séparé, $\pi : \mathfrak X\rightarrow M$ le morphisme vers son espace des modules grossier,
qu'on suppose projectif sur $S$, de dimension $n$, et muni d'un faisceau inversible très ample $\mathcal O(1)$. On note alors, pour $\mathcal F$ un faisceau localement libre sur $\mathfrak X$ :
$$\deg_{\mathfrak X}(\mathcal F)=\int_{\mathfrak X}^{et} c_1^{et}(\mathcal F)\cdot \pi^* c_1^{et}(\mathcal O(1))^{n-1}$$
\end{defi}

Pour des détails sur la théorie de l'intersection utilisée ici,
on renvoie à \ref{homcohom}. Intuitivement on peut comprendre ceci de la manière suivante : $\pi^* : \Pic M\otimes_{\mathbb Z} \mathbb Q\rightarrow \Pic \mathfrak X\otimes_{\mathbb Z} \mathbb Q$ est un isomorphisme, et l'on peut voir $\deg_{\mathfrak X}$ défini grâce à $\deg_M$.

\subsubsection{Comparaison des deux notions}

\begin{thm}
\label{thmdeg}
$r\geq 1$ un entier naturel, $X$ une variété projective lisse sur $k$ algébriquement clos de caractéristique première à $r$, $\mathcal O(1)$ un faisceau inversible très ample sur $X$, $D$ un diviseur de Cartier effectif sur $X$, lisse sur $k$, $\mathcal L=\mathcal O_X(D)$, $s$ la section canonique, $\mathfrak X=X_{(\mathcal L,s,r)}$ le champ des racines $r$-ièmes,  $\mathcal E_\cdot$ un objet de $\Par_{\frac{1}{r}} (X,D)$, $\mathcal F=G(\mathcal E_\cdot)$ le faisceau localement libre sur $\mathfrak X$ associé via la correspondance du théorème \ref{thmprinc}. Alors on a :

$$ \deg_{par}(\mathcal E_\cdot)=\deg_{\mathfrak X}(\mathcal F)$$

\end{thm}

\begin{rem}
On peut écrire formellement $\deg_{par}(\mathcal E_\cdot)=\int_{-1}^0 \deg(\mathcal E_t) dt$, ce qui a un sens, en définissant convenablement la mesure $dt$, ou bien en redéfinissant $\mathcal E_\cdot$ comme une fonction localement constante. On retrouve alors la définition \ref{degpar}. Malgré l'analogie des notations, je n'ai pas trouvé de démonstration simple du théorème \ref{thmdeg} en inversant les signes $\int$.
\end{rem}

\begin{cor}
Soient $\mathcal E_\cdot$,  $\mathcal E'_\cdot$ deux objets de
$\Par_{\frac{1}{r}} (X,D)$. On a :

$$\deg_{par}(\mathcal E_\cdot\otimes\mathcal E'_\cdot)=
\rk\mathcal E_\cdot \deg_{par}(\mathcal E'_\cdot)+\rk\mathcal E'_\cdot
 \deg_{par}(\mathcal E_\cdot)$$
\end{cor}

\begin{rem}
Lorsque $\dim X=1$, ce corollaire est énoncé comme une facile (?)
conséquence de la définition du produit tensoriel des fibrés paraboliques
dans \cite{Biswas}.
\end{rem}

\subsection{Preuve}

\subsubsection{Expression champêtre de  $\chi_{par}(\mathcal E_\cdot)$}

\begin{lem}
\label{intcarpar}
Avec les notations du théorème \ref{thmdeg}, si de plus $\mathcal N$ est la racine $r$-ième canonique de $\mathcal L$ sur le champ $\mathfrak X=X_{(\mathcal L,s,r)}$, on a alors :
$$\chi_{par}(\mathcal E_\cdot)=\frac{1}{r}\chi(\mathfrak X,\mathcal F  \otimes
\oplus_{l=1}^r\mathcal N^{\otimes -l})$$
\end{lem}

\begin{proof}
Il suffit d'observer que si $F: \Vect (X_{(\mathcal L,s,r)})\rightarrow\Par_{\frac{1}{r}} (X,D)$ désigne la correspondance du théorème \ref{thmprinc}, alors on a $F(\mathcal F\otimes \mathcal N^{\otimes -l})=F(\mathcal F)[\frac{l}{r}]$.
\end{proof}

\begin{rem}
On suppose que l'on est dans la situation du corollaire \ref{quotunif}, i.e. que l'on dispose d'un morphisme $p:Y\rightarrow \mathfrak X$, où $Y$ est un schéma, tel que le morphisme composé $\pi\circ p:Y\rightarrow X$ est un revêtement cyclique uniforme, induisant un isomorphisme $[Y| {\boldsymbol{\mu}_r}]\simeq \mathfrak X$. Alors la formule de projection et le lemme \ref{intcarpar}
donnent $\chi_{par}(\mathcal E_\cdot)
=\frac{1}{r}\chi(Y,p^*(\mathcal F  \otimes \mathcal N^{\otimes -1}))$, ce dont on conclut $ \deg_{par}(\mathcal E_\cdot)=\frac{1}{r}\deg_{Y}(p^*\mathcal F)$ (où
$\deg_{Y}(\cdot)$ est pris relativement à l'image réciproque sur $Y$ du fibré très ample $\mathcal O(1)$ sur $X$), d'où le théorème \ref{thmdeg} dans ce cas. Toutefois, en l'absence d'une telle présentation de $\mathfrak X$, on est contraint de procéder autrement.
\end{rem}

\subsubsection{Application de GRR}

Le théorème de Grothendieck-Riemann-Roch pour les champs de Deligne-Mumford, du à B.Toën (\cite{Toen}, voir aussi appendice \ref{GRR}),
permet d'obtenir une première expression pour $\deg_{par}(\mathcal E_\cdot)$.

Précisons tout d'abord une notation : l'isomorphisme de schémas 
$I_X\simeq X$ (voir \ref{inertie}) induit en cohomologie un isomorphisme d'algèbres $H_{et}^*(X)\simeq H_{et}^*(I_X)=H_{rep}^*(X)$, et pour $x\in K_0(X)$, on notera $c_i^{rep}(x)$ l'image de $c_i^{et}(x)$ par cet isomorphisme (c'est seulement dans ce cas que l'on s'autorisera à parler de classes de Chern à coefficients dans les représentations).

\begin{lem}
$$\deg_{par}(\mathcal E_\cdot)=\frac{1}{r}
\int_{\mathfrak X}^{rep}\td^{rep}(\mathfrak X)\cdot(\ch^{rep}(\mathcal F)-
\ch^{rep}(\mathcal O_{\mathfrak X}^{\oplus \rho}))\cdot\ch^{rep}(\oplus_{l=1}^r\mathcal N^{\otimes -l})\cdot\pi^* c_1^{rep}(\mathcal O(1))^{n-1} $$
\end{lem}

\begin{proof}
On note tout d'abord que puisque $\mathcal F=G(\mathcal E_{\cdot})$, 
on a $G(\mathcal E_{\cdot}(m))\simeq \mathcal F\otimes \pi^*\mathcal O(m)$.
D'après les corollaires \ref{pipropre}, \ref{pilisse}, le morphisme 
$\mathfrak X \rightarrow \Spec k$ est propre et lisse.
 La formule donnée est donc une conséquence directe de la définition \ref{degpar}, du lemme \ref{intcarpar}, du corollaire \ref{HRRDM} de l'appendice \ref{GRR},
et de l'expression standard
$$\ch^{et}(\mathcal O(m))=\sum_{l=0}^{n}c_1^{et}(\mathcal O(1))^l\frac{m^l}{l!}$$
\end{proof}

Posons pour simplifier

$$x=\td^{rep}(\mathfrak
X)\cdot(\ch^{rep}(\mathcal F)- \ch^{rep}(\mathcal O_{\mathfrak X}^{\oplus
  \rho}))\cdot\ch^{rep}(\oplus_{l=1}^r\mathcal N^{\otimes -l})\cdot
\pi^*c_1^{rep}(\mathcal O(1))^{n-1}$$

On sait d'après \ref{competrep} que :

$$\int_{\mathfrak X}^{rep} x=\int_{I_{\mathfrak X}}^{et}x=
\int_{\mathfrak X}^{et}x_1+\int_{I_{\mathfrak X}-\mathfrak X}^{et}x_{\neq 1}$$

Le premier terme se laisse facilement calculer :

\begin{lem}
$$\int_{\mathfrak X}^{et}x_1=r\deg_{\mathfrak X}\mathcal F$$
\end{lem}

\begin{proof}
Rappelons que $y\rightarrow y_1$ est un morphisme d'anneaux.
Le lemme résulte donc des lemmes \ref{chrepchet} et \ref{tdreptdet} de l'appendice \ref{GRR}, ainsi que du fait facile 
$(\pi^*c_1^{rep}(\mathcal O(1))^{n-1})_1=\pi^*c_1^{et}(\mathcal O(1))^{n-1}$.
\end{proof}

Le théorème \ref{thmdeg} résultera donc de la preuve de $\int_{I_{\mathfrak X}-\mathfrak X}^{et}x_{\neq 1}=0$.

\subsubsection{Étude du champ d'inertie}

Pour analyser l'expression ci-dessus on a besoin de comprendre la géométrie du champ d'inertie $I_{\mathfrak X}$ (voir \ref{inertie}) du champ $\mathfrak X=X_{(\mathcal L,s,r)}$. Pour ceci, on peut travailler dans un cadre plus général.

On fixe un entier $r$, et comme schéma de base $S_0=\Spec (\mathbb Z[r^{-1}](\boldsymbol{\mu}_r))$.
Donné un schéma noethérien $X\rightarrow S$, $\mathcal L$ un faisceau inversible sur $X$, on peut lui associer sa gerbe (sur $X$) des racines $r$-ièmes  de la manière suivante : on reprend la définition de
$X_{(\mathcal L,s,r)}$ en oubliant la condition sur les sections, c'est-à-dire :

\begin{defi}[\cite{Cadman}]

  On définit une catégorie $X_{(\mathcal L,r)}$ fibrée en groupoïdes
  sur $S$ par
\begin{enumerate}
\item  les objets de $X_{(\mathcal L,r)}$ sont les triplets
    $(f,\mathcal M,\phi)$ où $f:S\rightarrow X$ est un morphisme
    de schémas, $\mathcal M$ un faisceau inversible sur $S$,
et $\phi : \mathcal M^{\otimes r} \simeq f^*\mathcal L$ un isomorphisme,
\item les morphismes de $(f,\mathcal M,\phi)$ à $(g,\mathcal N,\psi)$ (tous deux au dessus de $S$) sont les isomorphisme de faisceaux $\nu : \mathcal
  M\simeq \mathcal N$ tels que $\psi \circ \nu^{\otimes r}=\phi$.
\end{enumerate}
  \end{defi}

On peut à présent énoncer une proposition à propos de  $I_{\mathfrak X}$ :

\begin{prop}
\label{strucinert}
Soit $r$ un entier, $X$ un schéma noethérien sur $S=\Spec (\mathbb Z[r^{-1}](\boldsymbol{\mu}_r))$, $\mathcal L$ un faisceau inversible sur $X$,
$s$ une section globale de $\mathcal L$, $D=Z(s)$ le lieu des zéros,
$\mathfrak X=X_{(\mathcal L,s,r)}$ le champ des racines $r$-ièmes de $\mathcal L$ sur $X$, $\mathfrak D=D_{(\mathcal L_{|D},r)}$ la gerbe
des racines $r$-ièmes de $\mathcal L_{|D}$ sur $D$.
On a alors un isomorphisme canonique :

$$I_{\mathfrak X} \simeq \mathfrak X \coprod \mathfrak D^{\amalg \boldsymbol{\mu}_r(S_0)-\{1\}}$$

\end{prop}

\begin{proof}

On note tout d'abord qu'on a un morphisme canonique
$\mathfrak i :\mathfrak D \rightarrow \mathfrak X$ composé des morphismes évidents
$D_{(\mathcal L_{|D},r)}\rightarrow D_{(\mathcal L_{|D},0,r)}$ et de
$D_{(\mathcal L_{|D},0,r)}\rightarrow X_{(\mathcal L,s,r)}$, comme les deux sont des immersions fermées ($D_{(\mathcal L_{|D},0,r)}$ est le $r$-ième voisinage infinitésimal de $D_{(\mathcal L_{|D},r)}$ dans $X_{(\mathcal L,s,r)}$, voir \cite{Cadman})
$\mathfrak D \rightarrow \mathfrak X$ l'est aussi.
Il en résulte un morphisme canonique $I_{\mathfrak D} \rightarrow I_{\mathfrak X}$, c'est d'ailleurs également une immersion fermée (car en fait obtenue par changement de base à partir de la précédente).

On note d'autre part qu'on a un morphisme canonique $I_{\mathfrak X}\rightarrow
\boldsymbol{\mu}_{r,S}$. En effet soit un objet $(f,\mathcal M,t,\phi)$ de $\mathfrak X$ au dessus de $f:S\rightarrow X$, alors l'isomorphisme canonique $a_\mathcal M :\HOM(\mathcal M,\mathcal M)\simeq \mathcal O_S$
induit un morphisme de schémas en groupes sur $S$ :
$$ \AUT_S(f,\mathcal M,t,\phi)\rightarrow \boldsymbol{\mu}_{r,S}$$

On en déduit le morphisme recherché, en associant à l'objet
$((f,\mathcal M,t,\phi),\nu)$ de $I_{\mathfrak X}$ au dessus de
$f:S\rightarrow X$ l'élément $a_\mathcal M(S)(\nu)$ de $\boldsymbol{\mu}_{r}(S)$.

Pour $\zeta : S_0\rightarrow \boldsymbol{\mu}_{r,S}$ on note
$I_{\mathfrak D}^{\zeta} \rightarrow I_{\mathfrak X}^{\zeta}$ la fibre correspondante
de $I_{\mathfrak D} \rightarrow I_{\mathfrak X}$. Vu le choix de la base $S_0$,
$\boldsymbol{\mu}_{r,S}$ est constant sur $S_0$, et donc $I_{\mathfrak X}=\coprod_{\zeta \in \boldsymbol{\mu}_{r}(S_0)}I_{\mathfrak X}^{\zeta}$.
Il est clair que $I_{\mathfrak X}^1$ est l'image de la section canonique $\mathfrak X\rightarrow I_{\mathfrak X}$. On vérifie d'autre part facilement que pour tout $\zeta$ dans  $\boldsymbol{\mu}_{r}(S_0)$, le composé
$I_{\mathfrak D}^{\zeta}\rightarrow I_{\mathfrak D} \rightarrow \mathfrak D$ est un isomorphisme (on a un inverse évident).
Il reste donc à vérifier que pour $\zeta\neq 1$ dans
$\boldsymbol{\mu}_{r}(S_0)$, le morphisme $I_{\mathfrak D}^{\zeta} \rightarrow I_{\mathfrak X}^{\zeta}$ est un isomorphisme. Au dessus de $S$, ce morphisme est donné par $((g,\mathcal M,\phi),\nu)\rightarrow ((i\circ g,\mathcal M,0,\phi),\nu)$
où $g :S\rightarrow D$ est un morphisme de schémas et $i:D\rightarrow X$ est l'inclusion canonique. Ceci définit clairement un foncteur fidèlement plein
$I_{\mathfrak D}^{\zeta}(S) \rightarrow I_{\mathfrak X}^{\zeta}(S)$. Pour vérifier qu'il est essentiellement surjectif, on fixe un objet $((f,\mathcal M,t,\phi),\nu)$ dans $I_{\mathfrak X}^{\zeta}(S)$. Donc $a_\mathcal M(S)(\nu)$ vaut $\zeta$ sur chaque composante connexe de $S$, et comme $1-\zeta$ est inversible sur $S_0$, $\id-\rho$ est inversible. On déduit alors de $\rho(t)=t$ le fait que $t=0$, et en conséquence $f^*s=t^{\otimes r}=0$, donc $f:S\rightarrow X$ se factorise à travers $i:D\rightarrow X$, comme souhaité.

\end{proof}

\subsubsection{Réduction}

Il résulte de la proposition \ref{strucinert} que l'immersion fermée
$\mathfrak i :\mathfrak D \rightarrow \mathfrak X$ induit un isomorphisme
$\tilde{\mathfrak i}=I\mathfrak i_{|I_{\mathfrak D}-\mathfrak D}:I_{\mathfrak D}-\mathfrak
D\rightarrow I_{\mathfrak X}-\mathfrak X$ et donc

$$\int_{I_{\mathfrak X}-\mathfrak X}^{et}x_{\neq 1}=
\int_{I_{\mathfrak D}-\mathfrak D}^{et}\tilde{\mathfrak i}^*(x_{\neq 1})=
\int_{I_{\mathfrak D}-\mathfrak D}^{et}({\mathfrak i}^*x)_{\neq 1}
$$

Or il est clair que $\int_{\mathfrak D}^{et}({\mathfrak i}^*x)_1=0$, en effet
ceci résulte du fait que $\dim \mathfrak D = n-1$, de 
$(\mathfrak i^*\pi^*c_1^{rep}(\mathcal O(1))^{n-1})_1=\mathfrak i^*\pi^*c_1^{et}(\mathcal O(1))^{n-1}$, et du lemme \ref{chrepchet}
qui entraîne $(\ch^{rep}(\mathfrak i^*\mathcal F)-
\ch^{rep}(\mathcal O_{\mathfrak D}^{\oplus \rho}))_1=
\ch^{et}(\mathfrak i^*\mathcal F)-\ch^{et}(\mathcal O_{\mathfrak D}^{\oplus
\rho})$, expression dont la composante de degré $0$ est nulle.
On conclut que $\int_{I_{\mathfrak X}-\mathfrak X}^{et}x_{\neq 1}=
\int_{\mathfrak D}^{rep}{\mathfrak i}^*x$, et il ne reste donc qu'à montrer la
nullité de cette dernière expression.

Cela résultera du lemme suivant. Remarquons que d'après preuve de 
la proposition
\ref{strucinert}, il existe un isomorphisme canonique $I_{\mathfrak D} \simeq
\mathfrak D^{\amalg \boldsymbol{\mu}_r(S_0)}$, qui induit un
isomorphisme $H^*_{rep}(\mathfrak D)\simeq H^*_{et}(\mathfrak D)^{\boldsymbol{\mu}_r(S_0)}$.

\begin{lem}
\label{lemcle}
Soit $\chi$ le caractère de la représentation canonique
$\boldsymbol{\mu}_r\subset \mathbb G_m$, et $\Lambda=\mathbb Q(\boldsymbol{\mu}_{\infty})$. Pour tout entier $l\in\mathbb Z$, l'isomorphisme canonique
$H^*_{rep}(\mathfrak D)_\Lambda\simeq H^*_{et}(\mathfrak D)_\Lambda^{\boldsymbol{\mu}_r(S_0)}$
envoie $\ch^{rep}(\mathfrak i^*\mathcal N^{\otimes l})$ sur
$\chi^l\ch^{et}(\mathfrak i^*\mathcal N^{\otimes l})$.
\end{lem}
\begin{proof}
Notons $\mathfrak i^*\mathcal N=\mathcal N_{|\mathfrak D}$, et $\pi_{\mathfrak
  D}:I_{\mathfrak D}\rightarrow {\mathfrak D}$ la projection canonique.
Pour calculer $\ch^{rep}(\mathcal N_{|\mathfrak D}^{\otimes l})$ on doit
essentiellement déterminer la décomposition de  $\pi_{\mathfrak
  D}^*\mathcal N_{|\mathfrak D}^{\otimes l}$ en sous-fibrés propres.
On note celle-ci $$\pi_{\mathfrak
  D}^*\mathcal N_{|\mathfrak D}^{\otimes l}
=\oplus_{\zeta'\in\boldsymbol{\mu}_{\infty}}(\pi_{\mathfrak D}^*\mathcal N_{|\mathfrak D}^{\otimes l})^{(\zeta')}$$

Le lemme alors résulte de la définition de $\ch^{rep}(\cdot)$ (et plus particulièrement de celle du caractère de Frobénius \ref{carfrob}) et du sous-lemme :
\begin{lem}
$$(\pi_{\mathfrak D}^*\mathcal N_{|\mathfrak D}^{\otimes
  l})^{(\zeta')}_{|I_{\mathfrak D}^{\zeta}}=
\left\{ \begin{array}{cr}
\mathcal N_{|I_{\mathfrak D}^{\zeta}}^{\otimes
  l}  & {\rm si}\;\; \zeta'=\zeta^l \\ \\
0 & {\rm sinon}
\end{array} \right .$$

\end{lem}

\begin{proof}
C'est à peu près tautologique vu que l'action de $\zeta$ (vu comme
automorphisme de $\mathfrak D$) sur $\mathcal N_{|I_{\mathfrak D}^{\zeta}}$ se
fait par multiplication par $\zeta$ (vu comme élément de $\Lambda$), et ceci d'après les définitions de $I_{\mathfrak D}$ et $\mathcal N$.
\end{proof}
\end{proof}

On peut achever la démonstration du théorème \ref{thmdeg}. 
En effet,
le lemme \ref{lemcle} implique que l'image de 
$\ch^{rep}(\oplus_{l=1}^r \mathfrak i^*\mathcal N^{\otimes
-l})$ par le morphisme d'anneaux composé 
$H^*_{rep}(\mathfrak D)_\Lambda\simeq H^*_{et}(\mathfrak D)_\Lambda^{\boldsymbol{\mu}_r(S_0)}\rightarrow H^0_{et}(\mathfrak D)_\Lambda^{\boldsymbol{\mu}_r(S_0)}$ est le caractère de la représentation régulière de $\boldsymbol{\mu}_r$. D'autre part, le lemme \ref{chrepchet} montre que l'image de $(\ch^{rep}(\mathfrak i^*\mathcal
F)-\ch^{rep}(\mathcal O_{\mathfrak D}^{\oplus \rho}))$ par ce même morphisme s'annule en $1$. On en conclut que le produit 
$\ch^{rep}(\oplus_{l=1}^r \mathfrak i^*\mathcal N^{\otimes
-l})\cdot (\ch^{rep}(\mathfrak i^*\mathcal
F)-\ch^{rep}(\mathcal O_{\mathfrak D}^{\oplus \rho}))$ s'envoie sur $0$ par 
$H^*_{rep}(\mathfrak D)_\Lambda\rightarrow H^0_{rep}(\mathfrak D)_\Lambda$, 
ce qui suffit à montrer
que $\int_{\mathfrak D}^{rep}{\mathfrak i}^*x=0$, comme souhaité.

\section{Fibrés paraboliques sur les courbes}

\label{courb}

Dans cette partie, on fixe $r\geq 1$ un entier, $X$ une courbe
projective et lisse sur un corps algébriquement clos $k$, $D$ un
diviseur effectif réduit sur $X$, $\mathcal L=\mathcal O_X(D)$, $s=s_D$
la section canonique, $\mathfrak X=X_{(\mathcal L,s,r)}$ le champ
des racines $r$-ièmes.

\subsection{Fibrés paraboliques semi-stables}

\begin{defi}
\label{sousfib}
Soit $\mathcal F$, $\mathcal F'$ deux faisceaux localement libres sur
$\mathfrak X$, et $\mathcal F'\rightarrow \mathcal F$ un
monomorphisme. On dit que $\mathcal F'$ est un \emph{sous-fibré} de
$\mathcal F$ si $\mathcal F/\mathcal F'$ est localement libre.
\end{defi}

\begin{rem}
\label{fibquot}
Pour conserver un vocabulaire symétrique (alors que la situation ne
l'est pas!) on parlera de \emph{fibré quotient} pour tout faisceau
localement libre $\mathcal F''$ quotient de $\mathcal F$.
\end {rem}

\begin{defi}
Soit $\mathcal F$ un faisceau localement libre sur $\mathfrak X$. On
définit sa \emph{pente} par $\mu(\mathcal F)=\deg_\mathfrak
X\mathcal F/\rk\mathcal F$. On dit que $\mathcal F$ est
\emph{semi-stable} s'il vérifie l'une des deux conditions équivalentes
suivantes :
\begin{enumerate}
\item pour tout sous-fibré $\mathcal F'$ de $\mathcal F$, on a  $\mu(\mathcal F')\leq \mu(\mathcal F)$,
\item pour tout fibré quotient $\mathcal F''$ de $\mathcal F$, on a  $\mu(\mathcal F)\leq \mu(\mathcal F'')$.
\end{enumerate}
\end{defi}

\begin{defi}
\label{fibeng} Soient $\mathcal F$, $\mathcal F'$ deux faisceaux
localement libres sur $\mathfrak X$, et $i :\mathcal F'\rightarrow
\mathcal F$ un monomorphisme. On désigne par $\mathcal
F'^{\hookrightarrow \mathcal F}$ le plus petit sous-fibré de
$\mathcal F$ contenant $\mathcal F'$.
\end{defi}

On peut définir $\mathcal F'^{\hookrightarrow \mathcal F}$ de
manière équivalente par $\mathcal F'^{\hookrightarrow \mathcal
F}=\ker(\mathcal F \twoheadrightarrow \coker i/\mathcal T (\coker
i)) $, où $\mathcal T(\cdot)$ désigne le sous-faisceau de torsion.
Donc l'injection $\mathcal F' \rightarrow \mathcal
F'^{\hookrightarrow \mathcal F}$ est génériquement un isomorphisme, et en particulier
 $\rk \mathcal F'^{\hookrightarrow \mathcal F}=\rk \mathcal F'$.
De plus, si $\pi : \mathfrak X\rightarrow X$ désigne la projection
canonique, on a : $\pi_*(\mathcal F'^{\hookrightarrow \mathcal
F})=(\pi_*\mathcal F')^{\hookrightarrow \pi_*\mathcal F}$.

\begin{lem}
\label{lemdeg}
Avec les notations de la définition \ref{fibeng} on a
: $\mu( \mathcal F')\leq \mu(\mathcal F'^{\hookrightarrow \mathcal
F})$
\end{lem}
\begin{proof}
Cela résulte du lemme suivant :
\begin{lem}
\label{lemdegbis}
Soit $\mathcal K'\rightarrow\mathcal K$ un monomorphisme de
faisceaux inversibles sur $\mathfrak X$. Alors $ \deg_\mathfrak
X\mathcal K'\leq\deg_\mathfrak X\mathcal K$.
\end{lem}
\begin{proof}
On peut supposer que $\mathcal K'=\mathcal O_\mathfrak X$. Notons
$D=P_1+\cdots P_m$, où les $P_i$ sont les points de $X$ dans le
support de $D$, $\mathcal L_i=\mathcal O_X(P_i)$ et $s_i$ les
sections canoniques correspondantes. D'après le lemme \ref{CartDisj}, on a 
$\mathfrak X\simeq  X_{\mathcal L_1,s_1,r}\times_X\cdots\times_X  X_{\mathcal
  L_m,s_m,r}$. Soit $\mathcal N_i$ l'image
réciproque de la racine $r$-ième canonique de $\mathcal L_i$ par
le morphisme $\mathfrak X\rightarrow X_{\mathcal
L_i,s_i,r}$. On peut vérifier de manière élémentaire l'existence 
d'entiers $r_i$, $1\leq i\leq m$, tels que $\mathcal K\simeq
\pi^*\pi_* \mathcal K\otimes \otimes_{i=1}^m \mathcal N_i^{\otimes
r_i}$ (ou encore utiliser le théorème \ref{thmprinc}). 
D'après le lemme \ref{PartEnt}, chaque $r_i$ est un entier compris entre $0$
et $r-1$, uniquement déterminé par cette écriture. Le monomorphisme $\mathcal
O_\mathfrak X\rightarrow\mathcal K$ donne en projection un
monomorphisme $\mathcal O_X\rightarrow\pi_*\mathcal K$, et le théorème de
Riemann-Roch sur $X$ implique $\deg_X\pi_*\mathcal K\geq 0$. 
Donc $ \deg_\mathfrak X\mathcal
K=\deg_X\pi_*\mathcal K+\sum_{i=1}^m\frac{r_i}{r}\geq 0$.

\end{proof}

\begin{rem}
Pour mesurer le confort donné par les champs, on pourra comparer cette preuve à sa version parabolique dans \cite{Boden}, Lemma 3.7.
\end{rem}

\end{proof}

\begin{defi}
Soit $\mathcal F$ un faisceau localement libre sur $\mathfrak X$.
\begin{enumerate}
\item Soit $\mathcal E$ un faisceau localement libre sur $X$ et
$\mathcal E\rightarrow \pi_*\mathcal F$ un monomorphisme vers le
faisceau sous-jacent à $\mathcal F$. On appelle \emph{sous-fibré
induit par  $\mathcal E$} le sous-fibré de $\mathcal F$ :
$\pi^*\mathcal E^{\hookrightarrow \mathcal F}$.
\item Soit $\mathcal E'$ un faisceau localement libre sur $X$ qui est quotient
$\pi_*\mathcal F\twoheadrightarrow \mathcal E'$ du faisceau
sous-jacent à $\mathcal F$. On appelle \emph{fibré quotient induit
par $\mathcal E'$} et on note $\pi^* \mathcal E'_{\twoheadleftarrow
\mathcal F}$ le fibré quotient de $\mathcal F$ : $\coker(\pi^*\mathcal E^{\hookrightarrow \mathcal F}\rightarrow \mathcal F)$, où $\mathcal E=\ker(\pi_*\mathcal F\twoheadrightarrow \mathcal E')$.
\end{enumerate}
\end{defi}

\begin{rem}
\begin{enumerate}
\item
Il est clair que $\pi_*(\pi^*\mathcal E^{\hookrightarrow \mathcal F})=
\mathcal E^{\hookrightarrow \pi_*\mathcal F}$ et $\pi_*(\pi^*
\mathcal E'_{\twoheadleftarrow \mathcal F})=\mathcal E'$. Il est facile de se convaincre qu'on retrouve les notions de sous-fibré parabolique induit et fibré
parabolique quotient induit de \cite{Seshadri}, voir définitions 
\ref{sousfibparind} et \ref{fibparquotind} de l'appendice \ref{defSeshadri}.
\item
\label{remsscomp}
  Le théorème \ref{thmdeg} et la remarque \ref{degcomp} de l'appendice 
\ref{defSeshadri} montrent que $\mathcal F$ est semi-stable si et
  seulement si le faisceau parabolique $\mathcal E_{\cdot}$ sur $(X,D)$ associé 
est semi-stable au sens de Seshadri. 
\end{enumerate}
\end{rem}

\subsection{Fibrés paraboliques finis}

\begin{defi}[\cite{NoriRFG}]
Un faisceau localement libre $\mathcal F$ sur $\mathfrak X$ est dit
\emph{fini} s'il existe deux polynômes \emph{distincts} $P,Q\in
\mathbb N[X]$ tels que $P(\mathcal F)\simeq Q(\mathcal F)$.
\end{defi}

\begin{prop}
Tout faisceau localement libre fini sur $\mathfrak X$ est semi-stable de
degré $0$.
\end{prop}

\begin{proof}
La preuve est mot pour mot celle des courbes usuelles \cite{NoriRFG}, \S 3, grâce à l'aide du lemme \ref{lemdeg}.
\end{proof}

On suppose désormais que $X=\mathbb P^1$, la droite projective, et
que $D=P_1+\cdots +P_m$.

\begin{thm}
\label{strucfibparfindrproj}
Soit $\mathcal E_\cdot$ un fibré parabolique \emph{fini} sur $(\mathbb
P^1,D)$, où $D$ est un diviseur effectif réduit de degré $m$.
Supposons que les poids associés à $\mathcal E_\cdot$ ont des
dénominateurs premiers à la caractéristique du corps de base $k$ et
que le faisceau sous-jacent se décompose sous la forme

$$\mathcal E_0=\oplus_{j=1}^{\rk \mathcal E}\mathcal O(d_j)$$

Alors pour tout $j$ on a : $-m<d_j\leq 0$.

\end{thm}

\begin{proof}
Le principe de la preuve suivante m'a été communiqué par I.Biswas.

Soit $\mathcal F$ le faisceau sur $\mathfrak X$ correspondant à
$\mathcal E_\cdot$. Il est fini, donc semi-stable de degré $0$.
L'hypothèse est donc que $\pi_*\mathcal F\simeq \oplus_{j=1}^{\rk \mathcal E}\mathcal O(d_j)$. Fixons $j$. On considère d'une part $\pi^* \mathcal O(d_j)^{\hookrightarrow \mathcal F}$ le sous-fibré de $\mathcal F$ induit par $\mathcal O(d_j)$. La semi-stabilité donne $\mu(\pi^*\mathcal O(d_j)^{\hookrightarrow \mathcal F})\leq 0$, ce qui implique, au vu du lemme \ref{lemdeg}, 
que $d_j\leq 0$. Soit d'autre part
$\pi^*\mathcal O(d_j)_{\twoheadleftarrow \mathcal F}$ le fibré quotient induit par $\mathcal O(d_j)$.
Le fait que $\mathcal F$ soit semi-stable de degré $0$ implique que
$\mu(\pi^*\mathcal O(d_j)_{\twoheadleftarrow \mathcal F})\geq 0$, ce qui implique, en écrivant comme dans la démonstration du lemme \ref{lemdegbis} 
$\pi^*\mathcal O(d_j)_{\twoheadleftarrow \mathcal F}\simeq
\pi^*\mathcal O(d_j)\otimes \otimes_{i=1}^m \mathcal N_i^{\otimes r_i}$, que $d_j>-m$.
\end{proof}

\appendix

\section{Fibrés paraboliques  : la définition de Seshadri}
\label{defSeshadri}
Dans cette partie, $X$ est une courbe projective et lisse sur un corps
algébriquement clos $k$, $D$ un diviseur effectif réduit sur $X$.

\subsection{Définition locale}
 \begin{defi}[\cite{Seshadri}]
   \label{def1}
   Un faisceau \emph{parabolique} $(\mathcal E,F_*,\alpha_*)$ sur
   $(X,D)$ est la donnée d'un faisceau localement libre $\mathcal E$ sur $X$,
   pour tout point $P$ du support de $D$ d'un drapeau de la
   fibre résiduelle $E_P:={\mathcal
   E}_P\otimes_{\mathcal{O}_{X,P}}k(P)$ :

   $$E_P=F_1(E_P)\supset F_2(E_P)\supset \cdots \supset F_{n_P}(E_P)
   \supset F_{n_P+1}(E_P) =0$$

   et une suite de nombres
   rationnels (poids) $(\alpha_{P,i})_{1\leq i\leq n_P}$ vérifiant
   $$ 0\leq \alpha_{P,1}<\cdots< \alpha_{P,n_P}<1$$

 $n_P$ est la longueur de la filtration en $P$ et si
   $m_{P,i}=\dim_{k(P)} (F_{i}(E_P)/F_{i+1}(E_P))$ alors
   la suite $(m_{P,i})_{1\leq i\leq n_P}$ est la suite des
   multiplicités en $P$.

    \end{defi}

\begin{defi}[\cite{Seshadri}]
   Soient $(\mathcal E,F_*,\alpha_*)$ et $(\mathcal E',F'_*,\alpha'_*)$
   deux faisceaux paraboliques sur $(X,D)$. Un morphisme $\phi
   :(\mathcal E,F_*,\alpha_*)\rightarrow (\mathcal E',F'_*,\alpha'_*)$
   de fibrés paraboliques est un morphisme $\phi :\mathcal
   E\rightarrow \mathcal E'$ de faisceaux vérifiant :

   $$\forall P\in |D|\;\; \forall i,j \;\;\;
   \alpha_{i,P}>\alpha'_{j,P} \implies \phi(F_i(E_P))\subset
   F'_{j+1}(E'_P)$$

  \end{defi}

\subsection{Filtration associée}

Soit $(\mathcal E,F_*,\alpha_*)$ un fibré parabolique sur
   $(X,D)$.

Il est bien connu (et apparemment du à C.Simpson, voir \cite{BodYok}) 
qu'on peut lui associer canoniquement une filtration
descendante $\mathcal E_\cdot=(\mathcal E_\alpha)_{\alpha\in \mathbb
  Q}$ de la manière suivante.

Pour $P\in |D|$ on pose par convention $\alpha_{0,P}=\alpha_{n_P,P}-1$
et $\alpha_{n_P+1}=1$.
On pose alors pour $1\leq i\leq n_P$ et   $\alpha_{i-1,P}<\alpha\leq
\alpha_{i,P}$ : $\mathcal E_\alpha^P =\ker (\mathcal E \rightarrow
E_P/F_i(E_P))$, puis

$$\mathcal E_\alpha = \cap_{P\in |D|} \mathcal E_\alpha^P$$
puis on étend la définition à tout $\alpha$ dans $\mathbb Q$ en
imposant $\mathcal E_{\alpha +l}=\mathcal E_\alpha(-lD)$ pour $l$ entier.

\begin{rem}
  Soient $(\mathcal E,F_*,\alpha_*)$ et $(\mathcal E',F'_*,\alpha'_*)$
   deux fibrés paraboliques sur $(X,D)$,
   $\mathcal E_\cdot$ et $\mathcal E'_\cdot$ les filtrations
   associées. Un morphisme de faisceaux $\phi :\mathcal
   E\rightarrow \mathcal E'$  définit un morphisme de faisceaux
   paraboliques si et seulement si
 $\forall \alpha  \in \mathbb Q \;\;\phi (\mathcal E_\alpha)\subset
   \mathcal E'_\alpha$.
\end{rem}

\subsection{Sous-faisceau parabolique et faisceau parabolique quotient induits}
On donne la version de Maruyama-Yokogawa, en termes des filtrations associées. 

\begin{defi}[\cite{MY}]
\label{sousfibparind}
Soit $\mathcal E_\cdot$ un faisceau parabolique sur $(X,D)$ et $\mathcal E'$ un sous-fibré de $\mathcal E_0$ (i.e. $\mathcal E'$ est un sous-faisceau de $\mathcal E_0$ et $\mathcal E_0/\mathcal E'$ est localement libre). Le \emph{sous-fibré parabolique induit} $\mathcal E'_\cdot$  est défini par 
$\mathcal E'_\cdot=\mathcal E'\cap\mathcal E_\cdot $.  
\end{defi}

\begin{defi}[\cite{MY}]
\label{fibparquotind}
Soit $\mathcal E_\cdot$ un faisceau parabolique sur $(X,D)$ et $g:\mathcal E_0\rightarrow \mathcal E''$ un épimorphisme vers le faisceau localement libre $\mathcal E''$. Le \emph{fibré parabolique quotient induit} $\mathcal E'_\cdot$ est défini par $\mathcal E''_\cdot=g(\mathcal E_\cdot) $.
\end{defi}

\subsection{Faisceau parabolique semi-stable}

\begin{defi}[\cite{Seshadri}]
 Soit $(\mathcal E,F_*,\alpha_*)$ un fibré parabolique sur $(X,D)$.
 On définit
 \begin{enumerate}
   \item son \emph{degré} par $\deg\mathcal E_\cdot=\deg \mathcal E+\sum_{P\in |D|}
  \sum_{i=1}^{n_P}m_{P,i}\alpha_{P,i} $

   \item sa \emph{pente} par $\mu(\mathcal E_\cdot) = \frac{\deg\mathcal
   E_\cdot}{\rk \mathcal E}$.
 \end{enumerate}
\end{defi}

\begin{rem}
\label{degcomp}
  Il est facile de vérifier que cette définition est compatible avec la
  définition \ref{degpar}. 
\end{rem}

\begin{defi}[\cite{Seshadri}]
Un fibré parabolique  $(\mathcal E,F_*,\alpha_*)$ sur $(X,D)$ est dit
\emph{semi-stable} s'il vérifie l'une des deux conditions équivalentes
suivantes :

\begin{enumerate}
\item Pour tout sous-fibré parabolique induit $\mathcal E'_\cdot$
  de $\mathcal E_\cdot$, on a $\mu(\mathcal E'_\cdot)\leq \mu(\mathcal E_\cdot)$
\item Pour tout fibré parabolique quotient induit $\mathcal E''_\cdot$
  de $\mathcal E_\cdot$, on a  $\mu(\mathcal E''_\cdot)\geq\mu(\mathcal E_\cdot)$.
\end{enumerate}

\end{defi}

\begin{rem}
Ce n'est pas exactement la définition de \cite{Seshadri}, mais cela lui est immédiatement équivalent, comme c'est d'ailleurs précisé dans la remarque 1. suivant la définition 6., p.69. La raison que nous avons d'éviter la première forme de la définition est qu'elle utilise une notion de sous-fibré parabolique (resp. fibré parabolique quotient) qui est n'est pas compatible (en fait plus générale) avec la définition champêtre précisée par la définition \ref{sousfib} (resp. par la remarque \ref{fibquot}).
\end{rem}

\section{Bouts}
\label{Bout}
On rappelle quelques notions tirées de (\cite{MacLane}).
\begin{defi}
  Soient $\mathcal I$, $\mathcal C$ deux catégories, $F:\mathcal
  I^{op}\times \mathcal I \rightarrow \mathcal C$ un foncteur de
  variance mixte.
\begin{enumerate}
  \item Un \emph{bout (wedge)} de $F$ à un objet $C$ de $\mathcal C$
  est une collection de flèches $\alpha_I : F(I,I)\rightarrow C$ dans
  $\mathcal C$ \emph{dinaturelle} au sens que pour tout flèche
  $f:I\rightarrow J$ dans $\mathcal I$, le diagramme suivant
  commute :
  \begin{center}
\xymatrix@R=2pt{
 &&& & F(I,I)\ar[dr]^{\alpha_I} & \\
 &&& F(J,I)\ar[ur]^{F(f,1)}\ar[dr]_{F(1,f)} && C \\
 &&& & F(J,J)\ar[ur]_{\alpha_J} & }
\end{center}

  \item Si un bout universel existe, on le note
 $\int^{I}F(I,I)$.
\end{enumerate}

\end{defi}

\begin{prop}[Fubini pour les bouts universels]

\label{Fubini}
Soient $\mathcal I$, $\mathcal J$ et $\mathcal C$ trois catégories,
$F:\mathcal I^{op}\times \mathcal I\times \mathcal J^{op}\times \mathcal J \rightarrow \mathcal C$ un
foncteur.
Si les bouts universels suivants existent, on a des isomorphismes
naturels :

$$\int^\mathcal I\int^\mathcal J F(I,I,J,J)\simeq \int^{\mathcal I \times\mathcal J}F(I,I,J,J)
\simeq \int^\mathcal J\int^\mathcal I F(I,I,J,J)$$
\end{prop}

\section{Grothendieck-Riemann-Roch pour les champs de Deligne-Mumford}
\label{GRR}

On résume les résultats de \cite{ToenThese}, \cite{Toen} que nous utilisons.
On note $\mathfrak X$ un champ de Deligne-Mumford noethérien.

\subsection{K-théorie}

La catégorie exacte $\Vect(\mathfrak X)$ permet de définir deux spectres de $K$-théorie, le spectre de $K$-théorie usuel $K(\mathfrak X)$ et le spectre de $K$-théorie étale $K_{et}(\mathfrak X)=H(\mathfrak X_{et},\underline K)$, où $H(\mathfrak X_{et},\cdot)$ désigne la cohomologie (généralisée) d'un préfaisceau de spectres $\cdot$ sur $\mathfrak X_{et}$ et $\underline K$ le préfaisceau de spectres de $K$-théorie sur $\mathfrak X_{et}$. Ils sont reliés par un morphisme canonique $\can : K(\mathfrak X)\rightarrow K_{et}(\mathfrak X)$.

On considérant à la place de $\Vect(\mathfrak X)$ la catégorie exacte $\Coh(\mathfrak X)$ des faisceaux cohérents sur $ \mathfrak X$ on obtient de manière analogue le spectre de $G$-théorie $G(\mathfrak X)$ et le spectre de $G$-théorie étale $G_{et}(\mathfrak X)$.

On dispose de morphismes de spectres canoniques $K(\mathfrak X)\rightarrow G(\mathfrak X)$ et $K_{et}(\mathfrak X)\rightarrow G_{et}(\mathfrak X)$.

\subsection{Théorie homologique-cohomologique}
\label{GRRet}
\subsubsection{Cohomologie}
\label{homcohom}
Soit $\mathcal K_i$ le faisceau associé au préfaisceau de groupes abéliens
$U\rightarrow  K_i(U)$ sur $\mathfrak X_{et}$.
On pose $$H^i_{et}(\mathfrak X)=H^i(\mathfrak X_{et},\mathcal K_i\otimes \mathbb Q)$$ On dispose alors d'une $\mathbb Q$-algèbre commutative graduée $H^*_{et}(\mathfrak X)$. $H^*_{et}(\mathfrak X)$ est contravariant en $\mathfrak X$.

La construction universelle de Gillet \cite{Gillet} donne des classes de Chern $c_{i}^{et}: K_{0,et}(\mathfrak X)\rightarrow H^i_{et}(\mathfrak X)$, et on définit un caractère de Chern $\ch^{et} : K_{0,et}(\mathfrak X)\rightarrow H^*_{et}(\mathfrak X)$ et une classe de Todd $\td^{et} : K_{0,et}(\mathfrak X)\rightarrow H^*_{et}(\mathfrak X)$ par les formules habituelles.

\subsubsection{Homologie}
On définit par ailleurs des groupes d'homologie de la manière suivante.

Pour un schéma $S$, et un entier $i$, soit :
$$\mathcal R^i(S) :\oplus_{x\in S^{(0)}}K_i(k(x))\rightarrow \oplus_{x\in S^{(1)}}K_{i-1}(k(x))\rightarrow \cdots\rightarrow \oplus_{x\in S^{(i)}}K_0(k(x))$$
le complexe de Gersten concentré en degrés $[-i,0]$. 
On note $\mathcal R^i$ le préfaisceau de complexes $S\rightarrow \mathcal R^i(S)$ sur 
$\mathfrak X_{et}$.
On pose : $H_i^{et}(\mathfrak X)= \mathbb H(\mathfrak X_{et},\mathcal 
R^i\otimes \mathbb Q)$.

On dispose de $H_*^{et}(\mathfrak X)$ qui est un $H^*_{et}(\mathfrak X)$-module gradué. $H_*^{et}(\mathfrak X)$ est covariant pour les morphismes propres.

Lorsque $\mathfrak X$ est lisse, il y a un isomorphisme naturel
$H^*_{et}(\mathfrak X)\simeq H_*^{et}(\mathfrak X)$, ce qui entraîne en
particulier que $H^i_{et}(\mathfrak X)=0$ pour $i>\dim \mathfrak X$.

Lorsque $p : \mathfrak X \rightarrow \Spec k$ est propre, on désigne $p_*$ par
$\int_{\mathfrak X}^{et}$ (noter que $H_*^{et}(\Spec k)\simeq \mathbb Q$).
Lorsque $\mathfrak X$ est de plus lisse, équidimensionel de dimension $n$, et muni d'un faisceau très ample $\mathcal O(1)$, ceci permet de définir le degré d'un faisceau localement libre $\mathcal F$ sur $X$ comme étant le nombre
rationnel : $$\deg_{\mathfrak X}(\mathcal F)=\int_{\mathfrak X}^{et} c_{1}^{et}(\mathcal F)\cdot c_{1}^{et}(\mathcal O(1))^{n-1}$$

Par ailleurs, on dispose d'un théorème de Grothendieck-Riemann-Roch pour la $G$-théorie \emph{étale} (\cite{ToenThese}, Théorème 3.33), qui peut s'énoncer ainsi : il existe une transformation naturelle (``caractère de Chern homologique étale'') $$\tau_{ \mathfrak X}^{et}:  G_{0,et}(\mathfrak X)\rightarrow H_*^{et}(\mathfrak X)$$ qui vaut $x\rightarrow \td^{et}(\mathcal T_{ \mathfrak X})\ch^{et}(x)$ si $\mathfrak X$ est lisse, où $\mathcal T_{ \mathfrak X}$ est le fibré tangent,
et qui est covariant pour les morphismes propres (non nécessairement représentables) de champs algébriques quasi-projectifs
$\mathfrak X\rightarrow\mathfrak X'$.

Ce théorème est essentiellement équivalent au théorème usuel pour les espaces de modules sous-jacents. La remarque clé due à B.Toën est que le morphisme 
$G_0(\mathfrak X)\rightarrow G_{0,et}(\mathfrak X)$ n'est pas covariant en général, mais qu'il est par contre possible de définir une transformation de type Riemann-Roch de $G_0(\mathfrak X)$ à $G_{0,et}(I_{\mathfrak X})$, où $I_{\mathfrak X}$ désigne le \emph{champ d'inertie} de $\mathfrak X$ : c'est ce que B.Toën nomme Lefschetz-Riemann-Roch, et que nous détaillons à présent.

\subsection{Lefschetz-Riemann-Roch}
\label{LRR}
\subsubsection{Champ d'inertie}
\label{inertie}
Le champ d'inertie $I_{\mathfrak X}$ de $\mathfrak X$ est défini par
$$I_{\mathfrak X}=\mathfrak X\times_ {\mathfrak X\times \mathfrak X}\mathfrak X$$

Moins formellement, il s'agit d'une catégorie fibrée en groupoïdes, dont les objets au-dessus d'un schéma $S$ sont les couples $(s,h)$, avec $s$ un objet de $\mathfrak X(S)$, et $h$ un élément de $\Hom_{\mathfrak X(S)}(s,s)$, et dont les flèches de $(s,h)$ vers $(s',h')$ ($s$, $s'$ objets de  $\mathfrak X(S)$) sont les morphismes $u\in \Hom_{\mathfrak X(S)}(s,s')$ tels que $u^{-1}h'u=h$.

$I_{\mathfrak X}$ est un champ en groupes sur $\mathfrak X$, et donc muni d'une section canonique $\mathfrak X\rightarrow I_{\mathfrak X}$ envoyant l'objet $s$
sur l'objet $(s,\id)$.

Lorsque $\mathfrak X=[Y|G]$ est un champ quotient, ce à quoi on peut se ramener localement, on a le diagramme commutatif suivant dont toutes les faces sont cartésiennes :

\xymatrix@R=2pt{
&G_Y\ar[ddd]\ar[ddr]\ar[rr]&&Y\ar[ddd]\ar[ddr]&\\
&&&&\\
&&G\times Y\ar[rr]\ar[ddd]&&Y\times Y\ar[ddd]\\
&I_{\mathfrak X}\ar[ddr]\ar[rr]&&\mathfrak X\ar[ddr]&\\
&&&& \\
&& \mathfrak X\ar[rr]&&\mathfrak X\times \mathfrak X}

Par définition, $G_Y$ est le stabilisateur de $Y\rightarrow Y$, un $Y$-groupe.
Ceci donne une présentation de $I_{\mathfrak X}=[G_Y|G]$ (noter que $G$ agit par conjugaison sur le premier facteur de $G_Y\subset G\times Y$, l'action préserve donc la structure de $Y$-groupe de $G_Y$, si bien que $I_\mathfrak X$ est un $\mathfrak X$-groupe). On peut encore détailler cette expression en décomposant $G_Y$ en ses composantes connexes.

\subsubsection{Caractère de Frobénius}
\label{carfrob}
Par la suite, on fait l'hypothèse simplificatrice suivante : on se donne un corps $k$ contenant les racines de l'unité et un champ $\mathfrak X \rightarrow \Spec k$ \emph{modéré}, i.e. l'ordre d'inertie de tout point de $ \mathfrak X$ est premier à la caractéristique de $k$.

Tout faisceau localement libre $\mathcal F$ sur $I_{\mathfrak X}$
se décompose canoniquement sous la forme $$\mathcal F=\oplus_{\zeta \in
\boldsymbol{\mu}_{\infty}}\mathcal F^{(\zeta)}$$
où  $\mathcal F^{(\zeta)}$ est le sous-fibré défini en restriction à $(s,h):
U\rightarrow\mathcal I_{\mathfrak X}$ de la manière suivante : l'action de $h$ sur $\mathcal F_{|U}$ est diagonalisable, et $\mathcal F^{(\zeta)}_{|U}$ est le sous-fibré propre pour la valeur propre $\zeta$.

En posant $\Lambda=\mathbb Q(\boldsymbol{\mu}_{\infty})$ et pour un groupe abélien $A$, $A_\Lambda=A\otimes_{\mathbb Z}\mathbb Q(\boldsymbol{\mu}_{\infty})$, on peut alors
définir un morphisme $\rho_{\mathfrak X}:K_0 (I_{\mathfrak X})_\Lambda\rightarrow
K_0 (I_{\mathfrak X})_\Lambda$ en posant $$\rho_{\mathfrak X}([\mathcal F])=
\sum_{\zeta \in \boldsymbol{\mu}_{\infty}}\zeta[\mathcal F^{(\zeta)}]$$

On définit alors le caractère de Frobénius $\phi_{\mathfrak X}$ de ${\mathfrak X}$ comme étant le composé

 \begin{center}
\xymatrix@R=2pt{
\phi_{\mathfrak X}:K_0 (\mathfrak X)_\Lambda\ar[r]^{\pi_{\mathfrak X}^*}&K_0 (I_{\mathfrak X})_\Lambda\ar[r]^{\rho_{\mathfrak X}}&K_0 (I_{\mathfrak X})_\Lambda\ar[r]^{can}&K_{0,et} (I_{\mathfrak X})_\Lambda
}
\end{center}
où $\pi_{\mathfrak X}^*$ est l'image réciproque via la projection
$\pi_{\mathfrak X} : I_{\mathfrak X} \rightarrow {\mathfrak X}$.

\subsubsection{Caractère de Frobénius homologique}

On définit ensuite l'analogue d'un caractère de Todd :
si $\Omega^1_{I_{\mathfrak X}/{\mathfrak X}}$ désigne le fibré conormal du morphisme non ramifié $\pi_{\mathfrak X} : I_{\mathfrak X} \rightarrow {\mathfrak X}$,
on pose $$\alpha_{\mathfrak X}=can\circ \rho_{\mathfrak X}(\lambda_{-1}([\Omega^1_{I_{\mathfrak X}/{\mathfrak X}}])) \in K_{0,et}(I_{\mathfrak X})_\Lambda$$

où comme d'habitude $\lambda_{-1}(x)=\sum_i(-1)^i\lambda_i(x)$. On vérifie que c'est un élément inversible de $K_{0,et}(I_{\mathfrak X})_\Lambda$.

Le théorème de Lefschetz-Riemann-Roch peut s'énoncer ainsi (voir \cite{ToenThese} Théorème 3.25 pour les détails). Il existe, pour $\mathfrak X \rightarrow \Spec k$ quasi-projectif, un caractère de Frobénius homologique
$\phi_{\mathfrak X}: G_0 (\mathfrak X)_\Lambda \rightarrow G_{0,et}(I_{\mathfrak X})_\Lambda$, qui vaut $x\rightarrow \alpha_{\mathfrak X}^{-1}\phi_{\mathfrak X}(x)$ lorsque
${\mathfrak X}$ est lisse, et qui est covariant pour les morphismes propres
 $\mathfrak X\rightarrow\mathfrak X'$ de dimension cohomologique finie lorsque $\mathfrak X$ est lisse.

\subsection{Grothendieck-Riemann-Roch}
\label{GRRT}
En combinant les deux transformations de Riemann-Roch de \ref{GRRet} et \ref{LRR} on est conduit aux définitions suivantes.

\subsubsection{Cohomologie et homologie à coefficients dans les représentations}

On pose donc logiquement $H^i_{rep}(\mathfrak X)=H^i_{et}(I_{\mathfrak X})$ et
$H_i^{rep}(\mathfrak X)=H_i^{et}(I_{\mathfrak X})$. Ces groupes héritent des propriétés de fonctorialité de leur analogues étales. En particulier lorsque $p:\mathfrak X\rightarrow \Spec k$ est propre, on définit $\int_{\mathfrak X}^{rep}=p_*$, comme d'habitude.

Le caractère de Chern $\ch^{rep}: K_0 (\mathfrak X)\rightarrow
H^*_{rep}(\mathfrak X)_\Lambda$ est défini par

$$\ch^{rep}(x)=\ch^{et}(\phi_{\mathfrak X}(x))$$

La classe de Todd est elle définie comme l'élément de
$H^*_{rep}(\mathfrak X)_\Lambda$ :

$$\td^{rep}(\mathfrak X)=\ch^{et}(\alpha_{\mathfrak X}^{-1})
\td^{et}(\mathcal T_{I_{\mathfrak X}})$$

\subsubsection{Caractère de Chern homologique à coefficients dans les représentations}

Dans cette section, on suppose de plus que l'espace des modules $M$ de $\mathfrak X$ est quasi-projectif, et ce pour tous les champs considérés.

Le théorème de Grothendieck-Riemann-Roch peut se résumer ainsi : il existe un
``caractère de Chern homologique à coefficients dans les représentations'' :
$$\tau_{ \mathfrak X}^{rep}:  G_{0}(\mathfrak X)_{\Lambda}\rightarrow H_*^{rep}(\mathfrak X)_{\Lambda}$$

qui vaut $x\rightarrow \td^{rep}(\mathfrak X)\ch^{rep}(x)$ si $\mathfrak X$ est lisse, qui est covariant pour les morphismes propres $\mathfrak X\rightarrow\mathfrak X'$.

En particulier on a :

\begin{cor}[Hirzebruch-Riemann-Roch pour les champs de Deligne-Mumford]
\label{HRRDM}
Soit $\mathcal F$ un faisceau cohérent sur un champ de Deligne-Mumford
$\mathfrak X\rightarrow \Spec k$ propre et lisse. Alors on a
$$\chi(\mathfrak X,\mathcal F)=\int_{\mathfrak X}^{rep}\td^{rep}(\mathfrak X)\ch^{rep}(\mathcal F)
$$
\end{cor}

\subsubsection{Lien entre cohomologie usuelle et cohomologie à coefficients
  dans les représentations}

\label{competrep}

Dans cette section, on énumère un certain nombre de faits techniques que nous
avons utilisés.

On note d'abord que les morphismes image réciproque le long de la section
unité $1:\mathfrak X\rightarrow I_{\mathfrak X}$ et de l'inclusion
$I_{\mathfrak X}-\mathfrak X\rightarrow  I_{\mathfrak X}$ induisent un
isomorphisme d'anneaux

$$H^*_{rep}(\mathfrak X)=H^*_{et}(I_{\mathfrak X})\simeq H^*_{et}(\mathfrak
X)\oplus H^*_{et}(I_{\mathfrak X}-\mathfrak X)$$

Pour $x\in H^*_{rep}(\mathfrak X)$ on notera $x=x_1+ x_{\neq 1}$ la
décomposition correspondante.

On a une décomposition analogue du groupe d'homologie
$H_*^{rep}(\mathfrak X)$, compatible avec la précédente via les isomorphismes
de Poincaré $H^*_{rep}(\mathfrak X)\simeq H_*^{rep}(\mathfrak X)$ lorsque $\mathfrak X$ est lisse.

Si $p:\mathfrak X\rightarrow \Spec k$ est propre, il en découle que :

 $$\int_{\mathfrak X}^{rep} x=\int_{I_{\mathfrak X}}^{et}x=
\int_{\mathfrak X}^{et}x_1+\int_{I_{\mathfrak X}-\mathfrak X}^{et}x_{\neq 1}$$

\begin{lem}
\label{chrepchet}
 Pour tout $y$ dans $K_0(\mathfrak X)_\Lambda$ et tout $i\geq 0$ on a :
$$\ch^{rep}(y)_1=\ch^{et}(y)$$
\end{lem}

\begin{proof}

Vu la définition des caractères de Chern, tout revient à montrer que
le diagramme suivant commute.

\xymatrix@R=2pt{
&&&K_0(\mathfrak X)_\Lambda\ar[r]^{\pi_{\mathfrak X}^*}\ar[ddr]_{id} & K_0(I_{\mathfrak X})_\Lambda \ar[r]^{\rho_{\mathfrak X}}&  K_0(I_{\mathfrak
  X})_\Lambda \ar[ldd]^{1^*}\ar[dd]^{can}\\
&&&&& \\
&&&&K_0(\mathfrak X)_\Lambda\ar[dd]^{can}&K_{0,et}(I_{\mathfrak X})_\Lambda\ar[dd]^{\ch^{et}} \\
&&&&& \\
&&&&K_{0,et}(\mathfrak X)_\Lambda\ar[dd]^{\ch^{et}}& H_{et}^*(I_{\mathfrak X})_\Lambda
\ar[ldd]^{1^*}\\
&&&&& \\
&&&&H_{et}^*(\mathfrak X)_\Lambda&
}

Il est clair que le parallélogramme du bas commute, quant au triangle du haut,
sa commutativité découle de la définition de $\rho_{\mathfrak X}$.

\end{proof}

\begin{lem}
\label{tdreptdet}
$$\td^{rep}(\mathfrak X)_1=\td^{et}(\mathfrak X)$$
\end{lem}

\begin{proof}
En effet $(\alpha_{\mathfrak X})_1=1^*can\circ \rho(\lambda_{-1}([\Omega^1_{I_{\mathfrak X}/{\mathfrak X}}]))=can\lambda_{-1}(0)=1$.
\end{proof}
\bibliography{ParRac}

\end{document}